\documentclass[11pt]{article}

\usepackage{amsmath}
\usepackage{amssymb}
\usepackage{latexsym}

\newtheorem{assumption}{Assumption}

\def\qed{ \ \vrule width.2cm height.2cm depth0cm\smallskip}
\def \ep{\hbox{ }\hfill$\Box$}

\newenvironment{proof}{\noindent {\bf Proof.\/}}{$\qed$\vskip 0.1in}

\newcommand{\la}{\langle}
\newcommand{\ra}{\rangle}

\newcommand{\eps}{\varepsilon}

\newcommand{\ba}{\begin{array}}
\newcommand{\ea}{\end{array}}
\newcommand{\be}{\begin{equation}}
\newcommand{\ee}{\end{equation}}
\newcommand{\bea}{\begin{eqnarray}}
\newcommand{\eea}{\end{eqnarray}}
\newcommand{\beaa}{\begin{eqnarray*}}
\newcommand{\eeaa}{\end{eqnarray*}}

\def\dbD{\mathbb{D}}
\def\dbE{\mathbb{E}}
\def\dbF{\mathbb{F}}

\def\dbH{\mathbb{H}}

\def\dbL{\mathbb{L}}

\def\dbP{\mathbb{P}}
\def\dbR{\mathbb{R}}
\def\dbS{\mathbb{S}}
\def\dbQ{\mathbb{Q}}

\def\a{\alpha}
\def\b{\beta}
\def\g{\gamma}
\def\d{\delta}
\def\e{\varepsilon}

\def\k{\kappa}
\def\l{\lambda}

\def\t{\tau}
\def\f{\varphi}
\def\th{\theta}
\def\o{\omega}

\def\G{\Gamma}

\def\L{\Lambda}

\def\O{\Omega}
\def\cF{{\cal F}}

\def\cL{{\cal L}}

\def\cP{{\cal P}}

\def\cY{{\cal Y}}
\def\cZ{{\cal Z}}
\def\no{\noindent}

\def\ms{\medskip}
\def\bs{\bigskip}
\def\q{\quad}

\def\pa{\partial}
\def\cd{\cdot}
\def\cds{\cdots}

\def\qed{ \hfill \vrule width.25cm height.25cm depth0cm\smallskip}

\newcommand{\basa}{\begin{assumption}}
\newcommand{\easa}{\end{assumption}}

\newcommand{\bas}{\begin{assum}}
\newcommand{\eas}{\end{assum}}

\def\limsup{\mathop{\overline{\rm lim}}}
\def\liminf{\mathop{\underline{\rm lim}}}

\def\esup{\mathop{\rm ess\;sup}}
\def\einf{\mathop{\rm ess\;inf}}

\def\pa{\partial}

 \def\cd{\cdot}
\def\cds{\cdots}

\def\dis{\displaystyle}

\def\cad{{c\`{a}dl\`{a}g}}

\def\1{\mathbf{1}}

\def\:{\!:\!}
\def\reff#1{{\rm(\ref{#1})}}
\def \proof{{\noindent \it Proof.\quad}}

at 9pt

\begin{document}

\newtheorem{thm}{Theorem}[section]
\newtheorem{lem}[thm]{Lemma}
\newtheorem{cor}[thm]{Corollary}
\newtheorem{prop}[thm]{Proposition}
\newtheorem{rem}[thm]{Remark}
\newtheorem{eg}[thm]{Example}
\newtheorem{defn}[thm]{Definition}
\newtheorem{assum}[thm]{Assumption}

\renewcommand {\theequation}{\arabic{section}.\arabic{equation}}
\def\thesection{\arabic{section}}

\title{Wellposedness of Second Order Backward SDEs}
\author{H.~Mete {\sc Soner}\footnote{ETH (Swiss Federal Institute of Technology),
Z\"urich and
Swiss Finance Institute, hmsoner@ethz.ch. Research partly supported by the
European Research Council under the grant 228053-FiRM.
Financial support from
 the ETH Foundation
is also gratefully acknowledged.}
        \and Nizar {\sc Touzi}\footnote{CMAP, Ecole Polytechnique Paris, nizar.touzi@polytechnique.edu. Research supported by the Chair {\it Financial Risks} of the {\it Risk Foundation} sponsored by Soci\'et\'e
               G\'en\'erale, the Chair {\it Derivatives of the Future} sponsored by the {F\'ed\'eration Bancaire Fran\c{c}aise}, and
               the Chair {\it Finance and Sustainable Development} sponsored by EDF and Calyon. }
        \and Jianfeng {\sc Zhang}\footnote{University of Southern California, Department of Mathematics, jianfenz@usc.edu. Research supported in part by NSF grant DMS 06-31366.}
\footnote{We are grateful to Marcel Nutz for his careful reading which helped clarifying some technical points of the proofs.}
}
\date{First version: April 5, 2010
      \\ This version: December 22, 2010}

\maketitle

\begin{abstract}
We provide an existence and uniqueness theory for an extension of backward SDEs to the second order. While standard Backward SDEs are naturally connected to semilinear PDEs, our second order extension is connected to fully nonlinear PDEs, as suggested in \cite{CSTV}. In particular, we provide a fully nonlinear extension of the Feynman-Kac formula. Unlike \cite{CSTV}, the alternative formulation of this paper insists that the equation must hold under a non-dominated family of mutually singular probability measures. The key argument is a stochastic representation, suggested by the optimal control interpretation, and analyzed in the accompanying paper \cite{STZ09c}.

\vspace{5mm}


\end{abstract}


\section{Introduction}
\setcounter{equation}{0}

Backward stochastic differential equations (BSDEs) appeared in Bismut \cite{B} in the linear case, and received considerable attention since the seminal paper of Pardoux and Peng \cite{PP}. The various developments are motivated by applications in probabilistic numerical methods for partial differential equations (PDEs), stochastic control, stochastic differential games, theoretical economics and financial mathematics.

On a filtered probability space $(\O,\cF,\{\cF_t\}_{t\in[0,1]},\dbP)$ generated by a Brownian motion $W$ with values in $\dbR^d$, a solution to a one-dimensional BSDE consists of a pair of progressively measurable processes $(Y,Z)$ taking values in $\dbR$ and $\dbR^d$, respectively, such that
 \beaa
 Y_t &=& \xi-\int_t^1 f_s(Y_s,Z_s) ds - \int_t^1 Z_s dW_s, \quad t \in [0,1],~~\dbP-\mbox{a.s.}
 \eeaa
where $f$ is a progressively measurable function from $[0,1] \times\Omega \times  \dbR \times \dbR^{d}$ to $\dbR$, and $\xi$ is an $\cF_1$-measurable random variable.

If the randomness in the parameters $f$ and $\xi$ is induced by the current value of a state process defined by a forward stochastic differential equation (SDE), then the BSDE is referred to as a Markov BSDE and its solution can be written as a deterministic function of time and the current value of the state process. For simplicity, we assume the forward process to be reduced to the Brownian motion, then under suitable regularity assumptions, this function can be shown to be the solution of a parabolic semilinear PDE.
 \beaa
 -\pa_tv - h^0(t,x,v,Dv,D^2v)=0
 &\mbox{where}&
 h^0(t,x,y,z,\g):=\frac12{\rm Tr}[\g]-f(t,x,y,z).
 \eeaa
In particular, this connection is the main ingredient for the Pardoux and Peng extension of the Feynman-Kac formula to semilinear PDEs. For a larger review of the theory of BSDEs, we refer to El Karoui, Peng and Quenez \cite{EPQ}.

Motivated by applications in financial mathematics and probabilistic numerical  methods for PDEs, Cheridito, Soner, Touzi and Victoir \cite{CSTV} introduced the notion of Second Order BSDEs (2BSDEs). The key issue is that, in the Markov case studied by \cite{CSTV}, 2BSDEs are connected to the larger class of fully nonlinear PDEs. This is achieved by introducing a further dependence of the generator $f$ on a process $\gamma$ which essentially identifies to the Hessian of the solution of the corresponding PDE. Then, a uniqueness result is proved in an appropriate set $\cZ$ for the process $Z$. The linear 2BSDE example reported in Section \ref{sect-linear 2BSDE} below shows clearly that the specification of the class $\cZ$ is crucial, and can not recover the natural class of square integrable processes, as in classical BSDEs. However, except for the trivial case where the PDE has a sufficiently smooth solution, the existence problem was left open in \cite{CSTV}.

In this paper, we provide a complete theory of existence and uniqueness for 2BSDEs. The key idea is a slightly different definition of 2BSDEs which consists in reinforcing the condition that the 2BSDE must hold $\dbP-$a.s. for every probability measure $\dbP$ in a non-dominated class of mutually singular measures introduced in Section \ref{sect-preliminary} below. The precise definition is reported in Section \ref{sect-2BSDE}. This new point of view is inspired from the quasi-sure analysis of Denis \& Martini \cite{DM} who established the connection between the so-called hedging problem in uncertain volatility models and the so-called Black-Scholes-Barrenblatt PDE. The latter is fully nonlinear and has a simple piecewise linear dependence on the second order term. We also observe an intimate connection between \cite{DM} and the $G-$stochastic integration theory of Peng \cite{Peng-G}, see Denis, Hu and Peng \cite{denishupeng}, and our paper \cite{STZ09b}.

In the present framework, uniqueness follows from a stochastic representation suggested by the optimal control interpretation. Our construction follows the idea of Peng \cite{Peng-G}. When the terminal random variable $\xi$ is in the space ${\rm UC}_b(\O)$ of bounded uniformly continuous maps of $\o$, the former stochastic representation is shown in our accompanying paper \cite{STZ09c} to be the solution of the 2BSDE . Then, we define the closure of ${\rm UC}_b(\O)$ under an appropriate norm. Our main result then shows that for any terminal random variable in this closure, the solution of the 2BSDE can be obtained as a limit of a sequence of solutions corresponding to bounded uniformly continuous final datum $(\xi_n)_n$. These are the main results of this paper and are reported in Section \ref{sect-wellposedness}.

Finally, we explore in Sections \ref{sect-PDE} and \ref{sect-dpp} the connection with fully nonlinear PDEs. In particular, we prove a fully nonlinear extension of the Feynman-Kac stochastic representation formula. Moreover, under some conditions, we show that the solution of a Markov 2BSDE is a deterministic function of the time and the current state which is a viscosity solution of the corresponding fully nonlinear PDE.



\section{Preliminaries}
\label{sect-preliminary}
\setcounter{equation}{0}

Let $\O:= \{\o\in C([0,1], \dbR^d): \o_0=0\}$ be the canonical space equipped with the uniform norm $\|\o\|_\infty := \sup_{0\le t\le 1}|\o_t|$, $B$ the canonical process, $\dbP_0$ the Wiener measure, $\dbF:= \{\cF_t\}_{0\le t\le 1}$ the filtration generated by $B$, and $\dbF^+:=\{\cF^+_t,0\le t\le 1\}$ the right limit of $\dbF$.

\subsection{The local martingale measures}

We say a probability measure $\dbP$ is a local martingale measure if
the canonical process $B$ is a local martingale under $\dbP$. By F\"{o}llmer \cite{follmer} (see also Karandikar \cite{Karandikar} for a more general result), there exists an $\dbF-$progressively measurable process, denoted as $\int_0^t B_s dB_s$, which coincides with the It\^o's integral, $\dbP-$a.s. for all local martingale measures $\dbP$. In particular, this provides a pathwise definition of
\beaa
\la B\ra_t := B_t B_t^{\rm T} - 2 \int_0^t B_s dB^{\rm T}_s &\mbox{and}& \hat a_t:= \limsup_{\e\downarrow 0} \frac{1}{\e}\Big(\la B\ra_t-\la B\ra_{t-\e}\Big),
\eeaa
where $^{\rm T}$ denotes the transposition, and the $\limsup$ is componentwise.
Clearly, $\la B\ra$ coincides with the $\dbP-$quadratic variation of $B$, $\dbP-$a.s. for all local martingale measures $\dbP$.

Let $\overline{\cP}_W$ denote the set of all local martingale measures $\dbP$ such that
\bea
\label{overlinecPW}
\la B\ra_t ~\mbox{is absolutely continuous in}~t
&\mbox{and}&
\hat a ~\mbox{takes values in}~ \dbS^{>0}_d, ~\dbP-\mbox{a.s.}
\eea
where $\dbS^{>0}_d$ denotes the space of all $d\times d$ real valued positive definite matrices. We note that, for different $\dbP_1, \dbP_2\in\overline{\cP}_W$, in general $\dbP_1$ and $\dbP_2$ are mutually singular. This is illustrated by the following example.

\begin{eg}
\label{example-singular}
{\rm Let $d=1$, $\dbP_1 := \dbP_0\circ (\sqrt{2}B)^{-1}$, and
$\O_i := \{\la B\ra_t = (1+i)t, t\ge 0\}$, $i=0, 1$. Then,
$\dbP_0, \dbP_1\in\overline{\cP}_W$,
$\dbP_0(\O_0) = \dbP_1(\O_1) = 1$, and $\dbP_0(\O_1)=\dbP_1(\O_0)=0$.
That is, $\dbP_0$ and $\dbP_1$ are mutually singular.
\ep}
\end{eg}

For any $\dbP\in \overline{\cP}_W$, it follows from the L\'evy characterization that the It\^o's stochastic integral under $\dbP$
\bea
\label{WP}
W^\dbP_t := \int_0^t \hat a^{-{1\slash 2}}_sdB_s, &t\in [0,1],&\dbP-\mbox{a.s.}
 \eea
defines a $\dbP-$Brownian motion.

This paper concentrates on the subclass $\overline\cP_S \subset \overline\cP_W$ consisting of all probability measures
 \bea
 \label{Xa}
 \dbP^\a := \dbP_0 \circ (X^\a)^{-1} &\mbox{where}& X^\a_t := \int_0^t \a_s^{1\slash 2} dB_s, t\in [0,1], \dbP_0-\mbox{a.s.}
 \eea
for some  $\dbF-$progressively measurable process $\a$ taking values in $\dbS^{>0}_d$ with
$\int_0^1 |\a_t| dt <\infty$, $\dbP_0-$a.s. With $\overline{\dbF}^\dbP$ (resp. $\overline{\dbF^{W^\dbP}}^\dbP$) denoting the $\dbP-$augmentation of the right-limit filtration generated by $B$ (resp. by $W^\dbP$), we recall from \cite{STZ09a} that
 \bea
 &\overline{\cP}_S = \big\{\dbP\in\overline{\cP}_W:  \overline{\dbF^{W^\dbP}}^\dbP = \overline{\dbF}^\dbP\big\},&
 \label{overlinecPS}
 \\
 &\mbox{and every}~\dbP\in\overline\cP_S~\mbox{satisfies}~\mbox{the Blumenthal zero-one law}&
 \nonumber
 \\
 &\mbox{and the martingale representation property.}&
 \label{0-1 MRP}
 \eea

\begin{rem}\label{rem-a1}{\rm
Let the process $\alpha$ be as above. Then by Lemma 2.2 in \cite{STZ09c},
\\
$\bullet$ there exists an $\dbF$-progressively measurable mapping $\b_\a$ such that $B_t = \b_\a(t,X^\a_.)$, $t\le 1$, $\dbP_0-$a.s.
\\
$\bullet$ the quadratic variation of the canonical process under $\dbP^\alpha$ is characterized by $\hat{a}(B)= \a\circ\b_\a(B)$, $dt\times\dbP^\a-$a.s.
}
\end{rem}

\begin{rem}\label{rem-a2}{\rm
As a consequence of the latter remark, given process $a$ with values in $\dbS_d^{>0}$ and $\int_0^1|a_t|dt<\infty$, it is not clear whether
there exists a process $\alpha$ as above 
so that the canonical process $\hat{a}=a$, $\dbP^\alpha-$a.s. The answer to this subtle question is negative in general, as shown by the example
 \beaa
 \alpha_t := \1_{[0,2]}(\hat{a}_t)
        +3~\1_{(2,\infty)}(\hat{a}_t),
 &t\in[0,1].&
 \eeaa
This will raise some technical problems in Section \ref{sect-PDE2}.
}
\end{rem}

\begin{rem}\label{rem-FF+}{\rm Let $\dbP\in\overline{\cP}_S$ be fixed. It follows from the Blumenthal zero-one law that $\dbE^\dbP[\xi|\cF_t]=\dbE^\dbP[\xi|\cF^+_t]$, $\dbP-$a.s. for any $t\in[0,1]$ and $\dbP-$integrable $\xi$. In particular, this shows that any $\cF_t^+-$measurable random variable has an $\cF_t-$measurable $\dbP-$mofication.
}
\end{rem}

\subsection{The nonlinear generator}

Our nonlinear generator is a map
 $$
 H_t(\o,y,z,\g):
 [0,1]\times \O\times \dbR\times \dbR^d\times D_H \to  \dbR,
 $$
where $D_H\subset\dbR^{d\times d}$ is a given subset containing $0$.
The corresponding conjugate of $H$ with respect to $\g$ takes values in $\dbR\cup\{\infty\}$ and is given by:
 \bea\label{F}
 \ba{c}
\dis F_t(\o,y,z,a)
 :=
 \sup_{\gamma\in D_H}~\left\{\frac{1}{2}a:\g - H_t(\o,y,z,\g)\right\}, ~~a\in \dbS^{>0}_d;\ms\\
\dis \hat F_t(y,z)
 :=
 F_t(y,z,\hat a_t)
\q \mbox{and}\q
 \hat F^0_t:=\hat F_t(0,0).
 \ea
 \eea
Here and in the sequel $a\:\g$ denotes the trace of the product matrix $a\g$.

We denote by $D_{F_t(y,z)}$ the domain of $F$ in $a$ for fixed $(t,\omega,y,z)$.

\begin{eg}
{\rm{The following are some examples of nonlinearities:
\\
1) Let $H_t(y,z,\gamma):=\frac12 a^0\:\gamma$ for some matrix $a^0\in\dbS_d^{>0}$. Here $D_H=\dbS_d$, and we directly calculate that $F_t(\o,y,z,a^0)=0$ and $F_t(y,z,a)=\infty$ whenever $a_t(\o)\neq a^0$. So $D_{F_t(y,z)}=\{a^0\}$.
\\
2) A more interesting nonlinearity considered by Peng \cite{Peng-G} will be commented later and is defined by $H_t(y,z,\gamma):=\frac12\sup_{a\in[\underline{a},\overline{a}]}(a\:\g)$. Here again $D_H=\dbS_d$, and we directly compute that $F_t(\o,y,z,a)=0$ for $a\in[\underline{a},\overline{a}]$, and $\infty$ otherwise. Hence $D_{F_t(y,z)}=[\underline{a},\overline{a}]$.
\\
3) Our last example is motivated by the problem of hedging under gamma constraints in financial mathematics. In the one-dimensional case, given two scalar $\underline{\Gamma}<0<\overline{\Gamma}$, the nonlinearity is $H_t(y,z,\gamma)=\frac12\gamma$ for $\gamma\in[\underline{\Gamma},\overline{\Gamma}]$, and $\infty$ otherwise. Here, $D_H=[\underline{\Gamma},\overline{\Gamma}]$ and $F_t(\o,y,z,a)=\frac12\big[\overline{\Gamma}(a-1)^+-\underline{\Gamma}(a-1)^-\big]$. In this example $D_{F_t(y,z)}=\dbR$.}}
\end{eg}

For the reason explained in Remark \ref{rem-norm} below, in this paper we shall fix a constant $\k$:
\bea
\label{k}
1 < \k \le 2,
\eea
and restrict the probability measures in the following subset  $\cP^\k_H\subset \overline{\cP}_S$:

\begin{defn}
 \label{defn-cP} Let $\cP^\k_H$ denote the collection of all those $\dbP\in\overline{\cP}_S$ such that
 \bea \label{ellipticity}
\underline a_\dbP \le \hat a \le \overline a_\dbP, ~~dt\times d\dbP-\mbox{a.s. for some}~ \underline a_\dbP, \overline a_\dbP\in \dbS^{>0}_d,~\mbox{and}~\dbE^\dbP\Big[\big(\int_0^1 |\hat F^0_t|^\k dt\big)^\frac{2}{\k} \Big] <\infty.
\eea
\end{defn}

It is clear that $\cP_H^\k$ is decreasing in $\k$, and $\hat a_t\in D_{F_t(0,0)}$, $dt\times d\dbP-$a.s. for all $\dbP\in{\cP_H^\k}$. Also, we emphasize on the fact that the bounds $(\underline a_\dbP,\overline a_\dbP)$ are not uniform in $\dbP$. In fact this restriction on the set of measure is not essential. For instance, if the nonlinearity (and the terminal data introduced later on) are bounded, then the bound is not needed.

\begin{defn}
\label{defn-qs}
We say a property holds ${\cP_H^\k}-$quasi-surely (${\cP_H^\k}-$q.s. for short) if it holds $\dbP-$a.s. for all $\dbP\in{\cP_H^\k}$.
\end{defn}

Throughout this paper, the nonlinearity is assumed to satisfy the following conditions.

\begin{assum}\label{assum-H}
$\cP^\k_H$ is not empty, and the domain $D_{F_t(y,z)}=D_{F_t}$ is independent of $(\o,y,z)$. Moreover, in $D_{F_t}$,
$F$ is $\dbF-$progressively measurable, uniformly continuous in $\o$ under the uniform convergence norm, and
 \bea\label{FLip-yz}
 |\hat F_t(y,z)-\hat F_t(y',z')|
 &\le&
 C\left(|y-y'|+|\hat a^{1/2}(z-z')|\right),
 ~~\cP^\k_H-\mbox{q.s.}
 \eea
for all $t\in[0,1]$, $y,y'\in\dbR$, $z,z'\in\dbR^d$.
\end{assum}

Clearly, one can formulate conditions on $H$ which imply the above Assumption. We prefer to place our assumptions on $F$ directly because this function will be the main object for our subsequent analysis.

\subsection{The spaces and norms}

We now introduce the spaces and norms which will be needed for the formulation of the second order BSDEs. Notice that all subsequent notations extend to the case $\k=1$.

For $p\ge 1$, $L^{p,\k}_H$ denotes the space of all $\cF_1-$measurable scalar r.v. $\xi$ with
 \beaa
 \|\xi\|_{L^{p,\k}_H}^p
 &:=&
 \sup_{\dbP\in{\cP_H^\k}} \dbE^{\dbP}[|\xi|^p]
 \;<\;
 \infty;
 \eeaa
$\dbH^{p,\k}_H$ denotes the space of all $\dbF^+-$progressively measurable $\dbR^d-$valued processes $Z$ with
  \beaa
  \|Z\|_{\dbH^{p,\k}_H}^p
  &:=&
  \sup_{\dbP\in{\cP_H^\k}} \dbE^{\dbP}\Big[\big(\int_0^1 |\hat a_t^{1\slash 2}Z_t|^2dt\big)^{p\slash 2}\Big]
  \;<\; \infty;
  \eeaa
$\dbD^{p,\k}_H$ denotes the space of all $\dbF^+-$progressively measurable $\dbR-$valued processes $Y$ with
   \beaa
   \cP_H^\k-\mbox{q.s. {\cad} paths},
   &\mbox{and}&
   \dis\|Y\|_{\dbD^{p,\k}_H}^p := \sup_{\dbP\in{\cP_H^\k}} \dbE^{\dbP}\Big[\sup_{0\le t\le 1}|Y_t|^p\Big]
   < \infty.
   \eeaa

For each $\xi\in L_H^{1,\k}$, $\dbP\in {\cP_H^\k}$, and $t\in [0,1]$, denote
 \beaa
 \dbE^{H, \dbP}_t[\xi]
 :=
 \esup_{\dbP' \in {\cP_H^\k}(t+, \dbP)}~^{\!\!\!\!\!\!\dbP}~
 \dbE^{\dbP'}_t[\xi]
 &\mbox{where}&
 {\cP_H^\k}(t+,\dbP)
 :=
 \{\dbP'\in {\cP_H^\k}: \dbP' = \dbP ~\mbox{on}~ \cF_t^+\}.
 \eeaa
It follows from Remark \ref{rem-FF+} that
 $\dbE^{\dbP}_t[\xi] := \dbE^\dbP[\xi|\cF_t] = \dbE^\dbP[\xi|\cF^+_t]$, $\dbP-$a.s.
  Then, for each $p\ge \k$,  we define
\begin{equation}
\label{L2*}
\dbL^{p,\k}_{H} := \big\{\xi\in L_H^{p,\k}: \|\xi\|_{\dbL^{p,\k}_{H}} <\infty\big\}~~\mbox{where}~~
\dis\|\xi\|^p_{\dbL^{p,\k}_{H}} := \sup_{\dbP\in{\cP_H^\k}} \dbE^\dbP\Big[\esup_{0\le t\le 1}~^{\!\!\!\!\!\!\dbP}~ \Big(\dbE^{H, \dbP}_t[|\xi|^\k]\Big)^{p/\k}\Big].
\end{equation}
The norm $\|\cd\|_{\dbL^{p,\k }_{H}}$ is somewhat less standard. Below, we justify this definition.

\begin{rem}
\label{rem-norm}
{\rm Assume  $\cP_H:=\cP^\k_H$ and $L^p_H:= L_H^{p,\k}$ do not depend on $\k$ (e.g. when $\hat F^0$ is bounded).
\\
(i) For $1\le \k_1 \le \k_2\le p$, it is clear that
 \beaa
 \|\xi\|_{L^{p}_H}\le \|\xi\|_{\dbL^{p,\k_1}_{H}}\le \|\xi\|_{\dbL^{p,\k_2}_{H}} &\mbox{and thus}& \dbL^{p,\k_2}_{H} \subset \dbL^{p,\k_1}_{H}
 \subset L^{p}_H.
 \eeaa
Moreover, as in our paper \cite{STZ09b} Lemma 6.2, under certain technical conditions, we have
 \beaa
 \|\xi\|_{\dbL^{p_1, p_1}_{H}} \le C_{p_2\slash p_1} \|\xi\|_{L^{p_2}_H} ~~\mbox{and thus}~~ L^{p_2}_H\subset \dbL^{p_1, p_1}_{H},
 &\mbox{for any}&
 1\le p_1 <p_2.
 \eeaa
(ii) In our paper \cite{STZ09b}, we used the norm $\|\cd\|_{\dbL^{p,1}_{H}}$. However, this norm does not work in the present paper due to the presence of the nonlinear generator, see Lemma \ref{lem-BSDEest}. So in this paper we shall assume $\k >1$ in order to obtain the norm estimates.

\ms

\no (iii) In the classical case where ${\cP_H}$ is reduced to a single measure ${\cP_H}=\{\dbP_0\}$, we have $\dbE^{H,\dbP_0}_t= \dbE^{\dbP_0}_t$ and the process $\{\dbE^{H,\dbP_0}_t[|\xi|^\k ],t\in[0,1]\}$ is a $\dbP_0-$martingale, then it follows immediately from the Doob's maximal inequality that, for all $1\le \k <p$,
\bea
\label{norm-equivalence}
\|\xi\|_{L^p(\dbP_0)}= \|\xi\|_{L^{p}_H} \le \|\xi\|_{\dbL^{p,\k }_H} \le C_{p,\k } \|\xi\|_{L^{p}_H} &\mbox{and thus}& \dbL_{H}^{p,\k }=L^{p}_H = L^p(\dbP_0).
\eea
However, the above equivalence does not hold when $\k =p$.
\ep}
\end{rem}

\begin{rem}
{\rm As in \cite{STZ09b}, in order to estimate $\|Y\|_{\dbD^{p,\k}_H}$ for the solution $Y$ to the 2BSDE with terminal condition $\xi$, it is natural to consider the supremum over $t$ in the norm of $\xi$. In fact we can show that the process $M_t:=\dbE^{H, \dbP}_t[|\xi|^\k ]$ a $\dbP-$supermartingale. Therefore it admits a {\cad} version and thus the term
$\sup_{t \in [0,1]} M_t$ is measurable.
\ep}
\end{rem}

Finally, we denote by ${\rm UC}_b(\O)$ the collection of all bounded and uniformly continuous maps $\xi:\O\longrightarrow\dbR$ with respect to the $\|.\|_\infty-$norm, and we let
\bea
\label{hatL2}
\cL_{H}^{p,\k } := \mbox{the closure of UC}_b(\O)~\mbox{under the norm}~ \|\cd\|_{\dbL^{p,\k }_{H}}, ~\mbox{for every} ~1\le \k \le p.
\eea
 Similar to  \reff{norm-equivalence}, we have
\begin{rem}
\label{rem-norm2}
{\rm In the  case ${\cP_H^\k} =\{\dbP_0\}$, we have $\cL^{p,\k }_H= \dbL_{H}^{p,\k }=L^{p,\k}_H = L^p(\dbP_0)$ for $1\le \k <p$.
}
\end{rem}

\section{The second order BSDEs}
\label{sect-2BSDE}
\setcounter{equation}{0}

We shall consider the following second order BSDE (2BSDE for short):
\bea
\label{2BSDE}
Y_t = \xi - \int_t^1 \hat F_s(Y_s, Z_s) ds - \int_t^1 Z_s dB_s + K_1 - K_t, ~~0\le t\le 1, ~~\mbox{${\cP_H^\k}-$q.s.}
\eea

\begin{defn}
\label{defn-2BSDE}
For $\xi \in\dbL^{2,\k}_{H}$, we say $(Y, Z)\in \dbD^{2,\k}_H \times \dbH^{2,\k}_H$ is a solution to 2BSDE \reff{2BSDE} if
\\
\noindent $\bullet$ $Y_T = \xi$, ${\cP_H^\k}-$q.s.
\\
\noindent $\bullet$ For each $\dbP\in {\cP_H^\k}$, the process $K^\dbP$ defined below has nondecreasing paths, $\dbP-$a.s.:
\bea
\label{KP}
K^\dbP_t := Y_0 - Y_t + \int_0^t \hat F_s(Y_s, Z_s) ds + \int_0^t Z_s dB_s, ~~0\le t\le 1,~~\dbP-\mbox{a.s.}
\eea
\noindent $\bullet$ The family $\{K^\dbP, \dbP\in{\cP_H^\k}\}$ defined in \reff{KP} satisfies the following minimum condition:
\bea
\label{minimum}
K^\dbP_t = \einf_{\dbP' \in {\cP_H^\k}(t+, \dbP)}~^{\!\!\!\!\!\!\!\!\!\!\!\dbP}~ \dbE^{\dbP'}_t[K^{\dbP'}_1],~~\dbP-\mbox{a.s. for all} ~\dbP\in {\cP_H^\k}, t\in [0,1].
\eea
Moreover, if the family $\{K^\dbP, \dbP\in{\cP_H^\k}\}$ can be aggregated into a universal process $K$, we call $(Y, Z, K)$ a solution of 2BSDE \reff{2BSDE}.
\end{defn}

Clearly, we may rewrite \reff{KP} as
\bea
\label{2BSDEP}
Y_t = \xi - \int_t^1 \hat F_s(Y_s, Z_s) ds - \int_t^1 Z_s dB_s + K^\dbP_1 - K^\dbP_t, ~~0\le t\le 1, ~~\dbP-\mbox{a.s.}
\eea
In particular, if $(Y, Z, K)$ is a solution of 2BSDE \reff{2BSDE} in the sense of the above definition, then it satisfies \reff{2BSDE} ${\cP_H^\k}-$q.s.

Finally, we note that, if $\dbP'\in {\cP_H^\k}(t+,\dbP)$, then $K^\dbP_s = K^{\dbP'}_s$, $0\le s\le t$, $\dbP-$a.s. and $\dbP'-$a.s.

\subsection{Connection with the second order stochastic target problem \cite{STZ09c}}

Let $(Y, Z)$ be a solution of 2BSDE \reff{2BSDE}. If the conjugate in \reff{F} has measurable maximizer, that is, there exists a process $\G$ such that
\bea
\label{Gamma}
\frac12 \hat a_t : \G_t - H_t (Y_t, Z_t, \G_t) = \hat F_t(Y_t, Z_t),
\eea
 then $(Y,Z, \G)$ satisfies
\bea
\label{2BSDE1}
Y_t = \xi - \int_t^1 \big[\frac12 \hat a_s : \G_s - H_s (Y_s, Z_s, \G_s)\big] ds - \int_t^1 Z_s dB_s + K_1 - K_t, 0\le t\le 1, \mbox{${\cP_H^\k}-$q.s.}
\eea
If $Z$ is a semi-martingale under each $\dbP\in\cP$ and $d \la Z, B\ra_t = \G_t d\la B\ra_t$, ${\cP_H^\k}-$q.s., then,
\bea
\label{2BSDE0}
Y_t = \xi + \int_t^1 H_s (Y_s, Z_s, \G_s)ds - \int_t^1 Z_s \circ dB_s + K_1 - K_t, ~~0\le t\le 1, ~~\mbox{${\cP_H^\k}-$q.s.}
\eea
Here $\circ$ denotes the Stratonovich integral. We note that \reff{2BSDE0}, \reff{2BSDE1}, and \reff{2BSDE} correspond to the second order target problem which was first introduced in \cite{ST} under a slightly different formulation. The present form, together with its first and second relaxations, were introduced in \cite{STZ09c}. In particular, in the Markovian case, the process $\G$ essentially corresponds to the second order derivative of the solution to a fully nonlinear PDE, see Section \ref{sect-PDE}. This justifies the denomination as "Second Order" BSDE of \cite{CSTV}. We choose to define 2BSDE in the form of \reff{2BSDE}, rather than \reff{2BSDE1} or \reff{2BSDE0}, because this formulation is most appropriate for establishing the wellposedness result, which is the main result of this paper and will be reported in Section \ref{sect-wellposedness} below.

\subsection{An alternative formulation of 2BSDEs}

In \cite{CSTV}, the authors investigate the following so called 2BSDE in Markovian framework:
\bea
\label{2BSDE-CSTV}
\left\{\ba{lll}
\dis Y_t = g(B_1) + \int_t^1 h(s, B_s, Y_s, Z_s, \G_s)ds - \int_t^1 Z_s \circ dB_s,\\
\dis d Z_t = \G_t dB_t + A_t dt,
\ea\right.
 ~~0\le t\le 1, ~~\dbP_0-\mbox{a.s.}
\eea
where $h$ is a deterministic function. Then uniqueness is proved in an appropriate space $\cZ$ for $Z$. The specification of $\cZ$ is crucial, and there can be no uniqueness result if the solution is allowed to be a general square integrable process. Indeed, the following "simplest" 2BSDE with $d=1$ has multiple solutions in the natural square integrable space:
\bea
\label{2BSDEc}
\left\{\ba{lll}
\dis Y_t = \int_t^1 \frac12 c\G_sds - \int_t^1 Z_s \circ dB_s,\\
\dis d Z_t = \G_t dB_t + A_t dt,
\ea\right.
 ~~0\le t\le 1, ~~\dbP_0-\mbox{a.s.}
\eea
where $c \neq 1$ is a constant. See Example \ref{eg-counterexample} below. The reason is that, unless $c=1$, $\dbP_0$ is not in ${\cP_H^\k}$ for $H(\g) := \frac12 c\g$. Also see subsection 3.4 below.

\subsection{Connection with $G-$expectations and $G-$martingales}

In \cite{STZ09b} we established the martingale representation theorem for $G-$martingales, which were introduced by Peng \cite{Peng-G}. In our framework, this corresponds to the specification $H_t(y,z,\gamma) = G(\g) := \frac12 \sup_{\underline a \le a \le \overline a} (a:\g)$, for some $\underline a, \overline a\in \dbS^{>0}_d$.

As an extension of \cite{STZ09b}, and as a special case of our current setting, we set
\bea
\label{H=G-f}
H_t(y,z,\g) := G(\g) - f_t(y,z).
\eea
Then one can easily check that:
\begin{itemize}
\item $D_{F_t} = [\underline a, \overline a]$ and $F_t(y,z,a) = f_t(y,z)$ for all $a\in [\underline a, \overline a]$;

\item{} $\dis{\cP_H^\k} = \left\{\dbP\in \overline\cP_s: \underline a\le \hat a\le \overline a, dt\times d\dbP-\mbox{a.s. and}~\dbE^\dbP\Big[\big(\int_0^1 |f_t(0,0)|^\k dt\big)^{\frac{2}{\k}}\Big]<\infty \right\}$.

\end{itemize}

In this case \reff{2BSDE} is reduced to the following 2BSDE:
\bea
\label{GBSDE}
Y_t &=& \xi + \int_t^1 f_s(Y_s, Z_s) ds - \int_t^1 Z_s dB_s + K_1 - K_t, ~~\mbox{${\cP_H^\k}-$q.s.}
\eea
Moreover, we may decompose $K$ into $d K_t = k_t dt + d K^0_t$, where $k\ge 0$ and $dK^0_t$ is a measure singular to the Lebesgue measure $dt$. One can easily check that there exists process $\G$ such that $G(\G_t)-\frac12\hat a_t :\G_t  = k_t$. Then \reff{GBSDE} becomes
\bea
\label{GBSDE1}
Y_t = \xi + \int_t^1\Big(\frac12\hat a_s :\G_s - G(\G_s) + f_s(Y_s, Z_s)\Big) ds - \int_t^1 Z_s dB_s + K^0_1 - K^0_t, ~\mbox{${\cP_H^\k}-$q.s.}
\eea
The wellposedness of the latter $G-$BSDE (with $K^0=0$ and $\k=2$) was left by Peng as an open problem.
We remark that, although the above two forms are equivalent, we prefer \reff{GBSDE} than \reff{GBSDE1} because the component $\G$ of the solution is not unique, and we have no appropriate norm for the process $\G$.

\subsection{Connection with the standard BSDE}

Let $H$ be the following linear function of $\g$:
\bea
\label{Hlinear}
H_t(y,z,\g) = \frac12 I_d\:\g - f_t(y,z),
\eea
where $I_d$ is is the identity matrix in $\dbR^d$. We remark that in this case we do not need to assume that $f$ is uniformly continuous in $\o$. Then, under obvious extension of notations, we have
 \beaa
 D_{F_t(\o)} = \{I_d\}
 &\mbox{and}&
 \hat F_t(y,z) = f_t(y,z).
 \eeaa
Assume that $\dbE^{\dbP_0}\big[\int_0^1 |f_t(0,0)|^2 dt\big]<\infty$, then ${\cP_H^\k} = {\cP_H^2}=\{\dbP_0\}$.
In this case, the minimum condition \reff{minimum} implies
\beaa
0 = K_0 = \dbE^{\dbP_0} [K_1] &\mbox{and thus}& K=0, ~~\dbP_0-\mbox{a.s.}
\eeaa
Hence, the 2BSDE \reff{2BSDE} is equivalent to the following standard BSDE:
\bea
\label{BSDE}
Y_t = \xi - \int_t^1 f_s(Y_s, Z_s) ds - \int_t^1 Z_s dB_s, ~~0\le t\le 1, ~~\dbP_0-\mbox{a.s.}
\eea
We note that, by Remark \ref{rem-norm2}, in this case we have
\beaa
\cL_{H}^{2,\k } = \dbL^{2,\k }_{H} = L^{2,\k}_H = \dbL^2(\dbP_0) &\mbox{for all}& 1\le \k <2.
\eeaa

\section{Wellposedness of 2BSDEs}
\label{sect-wellposedness}
\setcounter{equation}{0}

Throughout this paper Assumption \ref{assum-H} and the following assumption will always be in force.

\begin{assum}\label{assum-F}
The process $\hat F^0$ satisfies the integrability condition:
\bea
\label{F0-integrability}
\phi^{2,\k}_{H}
:=
\sup_{\dbP\in\cP^\k_H}
\dbE^\dbP\Big[ \esup_{0\le t\le 1}~^{\!\!\!\!\!\!\dbP}~ \big(\dbE^{H,\dbP}_t\big[\int_0^1|\hat F^0_s|^\k ds\big]\big)^{\frac{2}{\k}}\Big]
\;<\; \infty.
\eea
\end{assum}
Clearly the definition of $\phi^{2,\k}_{H}$ above is motivated by the norm $\|\xi\|_{\dbL^{2,\k}_H}$ in \reff{L2*}, and it satisfies
\bea
\label{F0-integrability2}
\sup_{\dbP\in\cP^\k_H}
\dbE^\dbP\Big[\Big(\int_0^1 |\hat F^0_t|dt\Big)^2\Big]
&\le&
\phi^{2,\k}_{H}.
\eea

For any $\dbP\in{\cP_H^\k}$, $\dbF^+-$stopping time $\t$, and $\cF^+_\t-$measurable random variable $\xi\in \dbL^2(\dbP)$, let $(\cY^\dbP, \cZ^\dbP) := (\cY^\dbP(\t,\xi), \cZ^\dbP(\t,\xi))$ denote the solution to the following standard BSDE:
\bea
\label{BSDEP}
\cY^\dbP_t = \xi - \int_t^\t \hat F_s(\cY^\dbP_s, \cZ^\dbP_s) ds - \int_t^\t \cZ^\dbP_s dB_s,~~0\le t\le \t, ~\dbP-\mbox{a.s.}
\eea
We have the following result which is slightly stronger than the standard ones in the literature. The proof is provided in subsection \ref{proof-lem-BSDEest} of the Appendix for completeness.

\begin{lem}
\label{lem-BSDEest}
Suppose  Assumption \ref{assum-H} holds.
Then,  for each $\dbP\in {\cP_H^\k}$, the BSDE \reff{BSDEP} has a unique solution satisfying
the following estimates:
\bea
\label{BSDEYest}
&& |\cY^\dbP_t|^2  \le C_\k \Big(\dbE^\dbP_t\Big[ |\xi|^\k + \int_t^1 |\hat F^0_s|^\k ds\Big]\Big)^{\frac{2}{\k}}, ~0\le t\le 1,~ \dbP-\mbox{a.s.}\\
\label{BSDEZest}
&& \dbE^\dbP\Big[\int_0^1 |\hat a_t^{1/2} \cZ^\dbP_t|^2 dt\Big]\le C_\k \dbE^\dbP\Big[\sup_{0\le t\le 1}\Big(\dbE^\dbP_t\Big[ |\xi|^\k + \int_0^1 |\hat F^0_s|^\k ds\Big]\Big)^{\frac{2}{\k}}\Big].
\eea
\end{lem}

We note that in above lemma, and in all subsequent results, we shall denote by $C$ a generic constant which may vary from line to line and depends only on the dimension $d$ and the Lipschitz constant in \reff{FLip-yz} of Assumption \ref{assum-H}. We shall also denote by $C_\k $ a generic constant which may depend on $\k $ as well. We emphasize that, due to the Lipschitz condition \reff{FLip-yz}, the constants $C$ and  $C_\k $ in the estimates will not depend on the bounds $\underline a_\dbP$ and $\overline a_\dbP$ in \reff{ellipticity}.

\subsection{Representation and uniqueness of the solution}

\begin{thm}
\label{thm-representation}
 Let Assumptions \ref{assum-H} and \ref{assum-F} hold. Assume that $\xi\in \dbL^{2,\k}_{H}$ and that $(Y, Z)\in \dbD^{2,\k}_H\times \dbH^{2,\k}_H$ is a solution to 2BSDE \reff{2BSDE}. Then, for any $\dbP\in{\cP_H^\k}$ and $0\le t_1  < t_2 \le 1$,
\bea
\label{DPP}
 Y_{t_1} = \esup_{\dbP' \in {\cP_H^\k}(t_1+, \dbP)}~^{\!\!\!\!\!\!\dbP}~ \cY^{\dbP'}_{t_1}(t_2,Y_{t_2}),~~\dbP-\mbox{a.s.}
\eea
Consequently, the 2BSDE \reff{2BSDE} has at most one solution in $\dbD^{2,\k}_H\times \dbH^{2,\k}_H$.
\end{thm}

\proof We first prove the last statement about uniqueness.
So suppose that \reff{DPP} holds.  Then as a special case with $t_2=1$ we obtain
\bea
\label{representation}
Y_t = \esup_{\dbP' \in {\cP_H^\k}(t+, \dbP)}~^{\!\!\!\!\!\!\dbP}~ \cY^{\dbP'}_{t}(1,\xi),~~\dbP-\mbox{a.s. for all} ~\dbP\in {\cP_H^\k}, t\in [0,1].
\eea
Therefore $Y$ is unique. To prove the uniqueness of $Z$, we observe that
 $d\la Y, B\ra_t = Z_t d\la B\ra_t$, ${\cP_H^\k}-$q.s..  Therefore the
 uniqueness of $Y$ implies that $Z$ is also unique.

It remains to prove \reff{DPP}.

\no (i)\quad Fix $0\le t_1 < t_2 \le 1$ and $\dbP\in {\cP_H^\k}$. For any $\dbP'\in {\cP_H^\k}(t_1+,\dbP)$, note that
\beaa
Y_t = Y_{t_2} - \int_t^{t_2}\hat F_s(Y_s, Z_s) ds - \int_t^1 Z_s dB_s + K^{\dbP'}_{t_2}-K^{\dbP'}_t,~~0\le t\le t_2, ~\dbP'-\mbox{a.s.}
\eeaa
and that $K^{\dbP'}$ is nondecreasing, $\dbP'-$a.s. By \reff{FLip-yz}, and applying the comparison principle for standard BSDE under $\dbP$, we have $Y_{t_1} \ge \cY^{\dbP'}_{t_1}(t_2, Y_{t_2})$, $\dbP'-$a.s. Since $\dbP'=\dbP$ on $\cF^+_{t_1}$, we get $Y_{t_1} \ge \cY^{\dbP'}_{t_1}(t_2, Y_{t_2})$, $\dbP-$a.s. and thus
\bea
\label{rep-est1}
Y_{t_1} \ge \esup_{\dbP' \in {\cP_H^\k}(t_1+, \dbP)}~^{\!\!\!\!\!\!\!\!\dbP}~ \cY^{\dbP'}_{t_1}(t_2, Y_{t_2}),~~\dbP-\mbox{a.s.}
\eea
\no (ii)\quad We now prove the other direction of the inequality. Fix $\dbP\in{\cP_H^\k}$.
For every $\dbP'\in{\cP_H^\k}(t_1+,\dbP)$, denote:
 \beaa
 \d Y := Y - \cY^{\dbP'}(t_2, Y_{t_2})
 &\mbox{and}&
 \d Z := Z - \cZ^{\dbP'}(t_2, Y_{t_2}).
 \eeaa
By the Lipschitz conditions \reff{FLip-yz}, there exist bounded processes $\l, \eta$ such that
\bea
\label{dYN}
\d Y_t = \int_t^{t_2} \big(\l_s \d Y_s + \eta_s \hat a_s^{1\slash 2} \d Z_s \big) ds
         - \int_t^{t_2} \d Z_s dB_s + K^{\dbP'}_{t_2}-K^{\dbP'}_t, ~t\le t_2, \dbP'-\mbox{a.s.}
\eea
Define:
\bea
\label{M}
M_t := \exp\Big(-\int_0^t \eta_s\hat a_s^{-1\slash 2} dB_s - \int_0^t ( \l_s + \frac12|\eta_s|^2) ds\Big),~~0\le t\le t_2, ~\dbP'-\mbox{a.s.}
\eea
By It\^o's formula, we have:
\bea
\label{MdY}
d\big(M_t \d Y_t\big)
=
M_t\big( \d Z_t - \d Y_t \eta_t{\hat a}^{-1/2}_t \big) dB_t
- M_t dK^{\dbP'}_t,
~~t_1\le t\le t_2, ~\dbP'-\mbox{a.s.}
\eea
Then, since $\d Y_{t_2}=0$, using standard localization arguments if necessary, we compute that:
\beaa
Y_{t_1}-\cY_{t_1}^{\dbP'}(t_2,Y_{t_2})
=
\d Y_{t_1}
=
M_{t_1}^{-1}\dbE^{\dbP'}_{t_1}\Big[\int_{t_1}^{t_2} M_t dK^{\dbP'}_t \Big]
\le
\dbE^{\dbP'}_{t_1}\Big[\sup_{t_1\le t\le t_2}(M_{t_1}^{-1} M_t) (K^{\dbP'}_{t_2}-K^{\dbP'}_{t_1})\Big]
\eeaa
by the non-decrease of $K^{\dbP'}$. By the boundedness of $\l, \eta$, for every $p\ge 1$ we have,
\bea
\label{Mest}
\dbE^{\dbP'}_{t_1}\Big[\sup_{t_1\le t\le t_2}(M_{t_1}^{-1} M_t)^p+\sup_{t_1\le t\le t_2}(M_{t_1} M_t^{-1})^p\Big] \le C_p,~~t_1\le t\le t_2,~\dbP'-\mbox{a.s.}
\eea
Then it follows from the H\"older inequality that:
 \beaa
 Y_{t_1}-\cY_t^{\dbP'}(t_2,Y_{t_2})
 &\le&
 \Big(\dbE^{\dbP'}_{t_1}\Big[\sup_{t_1\le t\le t_2} (M_{t_1}^{-1} M_t)^3
                        \Big]
 \Big)^{1/3}
 \Big(\dbE^{\dbP'}_{t_1}\big[(K^{\dbP'}_{t_2}-K^{\dbP'}_{t_1})^{3/2}
                        \big]
 \Big)^{2/3}
 \\
 &\le&
 C\Big(\dbE^{\dbP'}_{t_1}\big[K^{\dbP'}_{t_2}-K^{\dbP'}_{t_1}
                         \big]
       \dbE^{\dbP'}_{t_1}\big[(K^{\dbP'}_{t_2}-K^{\dbP'}_{t_1}\big)^2
                         \big]
 \Big)^{1/3}
 \eeaa
We shall prove in Step (iii) below that
 \bea\label{CPt1}
 C^\dbP_{t_1}:=
 \esup_{\dbP'\in{\cP_H^\k}(t_1+,\dbP)}~^{\!\!\!\!\!\!\!\!\dbP}~
 \dbE^{\dbP'}_{t_1}\left[(K^{\dbP'}_{t_2}-K^{\dbP'}_{t_1})^2\right]
 &<&
 \infty,~~\dbP-\mbox{a.s.}
 \eea
 Then, it follows from the last inequality that
 \beaa
 Y_{t_1}
 -\esup_{\dbP'\in\cP_H^\k(t_1+,\dbP)}\cY_{t_1}^{\dbP'}(t_2,Y_{t_2})
 &\le&
 C (C^\dbP_{t_1})^{1/3}
 \einf_{\dbP'\in\cP_H^\k(t_1+,\dbP)}
 \Big(\dbE^{\dbP'}_{t_1}\big[K^{\dbP'}_{t_2}-K^{\dbP'}_{t_1}\big]
 \Big)^{1/3}
 \;=\; 0,~~\dbP-\mbox{a.s.}
 \eeaa
by the minimum condition \reff{minimum}.

\no (iii)\quad  It remains to show that the estimate \reff{CPt1} holds. By the definition of the family $\{K^\dbP,\dbP\in{\cP_H^\k}\}$ we have:
\bea\label{est1-Neveu}
\sup_{\dbP'\in{\cP_H^\k}(t_1+,\dbP)}\dbE^{\dbP'}\left[(K^{\dbP'}_{t_2}-K^{\dbP'}_{t_1})^2\right]
&\le&
C\left(\|Y\|_{\dbD^{2,\k}_H}^2 + \|Z\|_{\dbH^{2,\k}_H}^2 + \phi^{2,\k}_{H}\right)<\infty.
\eea
We next use the definition of the essential supremum, see e.g. Neveu \cite{Neveu} to see that
 \bea\label{Neveu}
 \esup_{\dbP'\in{\cP_H^\k}(t_1+,\dbP)}
 \dbE^{\dbP'}_{t_1}\left[(K^{\dbP'}_{t_2}-K^{\dbP'}_{t_1})^2\right]
 &=&
 \sup_{n\ge 1}\dbE^{\dbP_n}_{t_1}\left[(K^{\dbP_n}_{t_2}-K^{\dbP_n}_{t_1})^2\right],~~\dbP-\mbox{a.s.}
 \eea
for some sequence $(\dbP_n)_{n\ge 1}\subset{\cP_H^\k}(t_1+,\dbP)$. Observe that for $\dbP'_1,\dbP'_2\in{\cP_H^\k}(t_1+,\dbP)$, there exists $\dbP'\in\cP^\k_H(t_1+,\dbP)$ such that
 \bea\label{filtrating}
 \dbE^{\dbP'}_{t_1}\left[(K^{\dbP'}_{t_2}-K^{\dbP'}_{t_1})^2\right]
 =
 \mu_{t_1}
 :=\max\left\{\dbE^{\dbP'_1}_{t_1}\left[(K^{\dbP'_1}_{t_2}
                                        -K^{\dbP'_1}_{t_1})^2\right],
       \dbE^{\dbP'_2}_{t_1}\left[(K^{\dbP'_2}_{t_2}
                                  -K^{\dbP'_2}_{t_1})^2\right]
     \right\}.
 \eea
Indeed, set
 \beaa
 E_1
 :=\Big\{\mu_{t_1}
 =\dbE^{\dbP'_1}_{t_1}\left[(K^{\dbP'_1}_{t_2}
                             -K^{\dbP'_1}_{t_1})^2\right]
   \Big\}
 &\mbox{and}&
 E_2:=\O\setminus E_1,
 \eeaa
 so that both sets are in $\cF_{t_1}$.  We then
 define the probability measure $\dbP'$ by,
 \beaa
 \dbP'[E] := \dbP'_1[E\cap E_1]+\dbP'_2[E\cap E_2]
 &\mbox{for all}&
 E\in\cF_1.
 \eeaa
Then, by its definition,  $\dbP'$ satisfies \reff{filtrating} trivially.
Moreover, in subsection \ref{proof-barPN-in-cP} of the Appendix, it is proved that
 \bea\label{barPN-in-cP}
 \dbP' &\in& {\cP_H^\k}(t_1+,\dbP) .
 \eea
Using this construction, by using a subsequence, if necessary,
we rewrite \reff{Neveu}, as
\beaa
 \esup_{\dbP'\in{\cP_H^\k}(t_1+,\dbP)}
 \dbE^{\dbP'}_{t_1}\left[(K^{\dbP'}_{t_2}-K^{\dbP'}_{t_1})^2\right]
 &=&
 \lim_{n\to\infty}\uparrow
 \dbE^{\dbP_n}_{t_1}\left[(K^{\dbP_n}_{t_2}-K^{\dbP_n}_{t_1})^2\right].
 \eeaa
It follows from \reff{est1-Neveu} that
 \beaa
 \dbE^\dbP\Big[\esup_{\dbP'\in{\cP_H^\k}(t_1+,\dbP)}
               \dbE^{\dbP'}_{t_1}\left[(K^{\dbP'}_{t_2}
                                        -K^{\dbP'}_{t_1})^2\right]
          \Big]
 &=&
 \dbE^\dbP\Big[\lim_{n\to\infty}\uparrow
               \dbE^{\dbP_n}_{t_1}\left[(K^{\dbP_n}_{t_2}
                                        -K^{\dbP_n}_{t_1})^2\right]
          \Big]
 \\
 &=&
 \lim_{n\to\infty}\uparrow
 \dbE^{\dbP_n}\left[(K^{\dbP_n}_{t_2}-K^{\dbP_n}_{t_1})^2\right]
 \\
 &\le&
 \sup_{\dbP'\in\cP^k_H(t_1+,\dbP)}
 \dbE^{\dbP'}\left[(K^{\dbP'}_{t_2}-K^{\dbP'}_{t_1})^2\right]
 \;<\;
 \infty
\eeaa
by \reff{est1-Neveu}, which implies the required estimate \reff{CPt1}.
\ep

\quad

As an immediate consequence of the representation formula \reff{representation}, together with the comparison principle for BSDEs, we have the following comparison principle for 2BSDEs.

\begin{cor}
\label{cor-comparison}
Let Assumptions \ref{assum-H} and \ref{assum-F} hold. Assume $\xi^i\in \dbL^{2,\k}_{H}$ and $(Y^i, Z^i)\in \dbD^{2,\k}_H\times \dbH^{2,\k}_H$ is a corresponding solution of the 2BSDE \reff{2BSDE}, $i=1,2$. If $\xi^1 \le \xi^2$, ${\cP_H^\k}-$q.s. then $Y^1 \le Y^2$, ${\cP_H^\k}-$q.s.
\end{cor}

\subsection{A priori estimates and the existence of a solution}

\begin{thm}
\label{thm-apriori}
Let Assumptions \ref{assum-H} and \ref{assum-F} hold.
\\
{\rm (i)}\quad Assume that $\xi\in \dbL^{2,\k}_{H}$ and that $(Y, Z)\in \dbD^{2,\k}_H\times \dbH^{2,\k}_H$ is a solution to 2BSDE \reff{2BSDE}. Then there exist a constant $C_\k$ such that
\bea
\label{est1}
 \|Y\|^2_{\dbD^{2,\k}_H} + \|Z\|^2_{\dbH^{2,\k}_H} + \sup_{\dbP\in{\cP_H^\k}} \dbE^\dbP[|K^\dbP_1|^2]
 &\le&
 C_\k\big(\|\xi\|^2_{\dbL^{2,\k}_{H}} + \phi^{2,\k}_{H}\big).
\eea
{\rm (ii)}\quad Assume that $\xi^i\in \dbL^{2,\k}_{H}$ and that $(Y^i, Z^i)\in \dbD^{2,\k}_H\times \dbH^{2,\k}_H$ is a corresponding solution to 2BSDE \reff{2BSDE}, $i=1,2$. Denote $\d \xi := \xi^1-\xi^2$, $\d Y := Y^1-Y^2$, $\d Z:= Z^1-Z^2$, and $\d K^\dbP := K^{1,\dbP}-K^{2,\dbP}$. Then there exists a constant $C_\k$ such that
\bea\label{est2}
\begin{array}{rcl}
 \dis\|\d Y\|_{\dbD^{2,\k}_H} &\le&  C_\k\|\d \xi\|_{\dbL^{2,\k}_{H}},
 \\
 \dis\|\d Z\|^2_{\dbH^{2,\k}_H} + \sup_{\dbP\in{\cP_H^\k}} \dbE^\dbP\Big[\sup_{0\le t\le 1}|\d K^\dbP_t|^2\Big]
 &\le&
 C_\k\|\d \xi\|_{\dbL^{2,\k}_{H}}\!\!\Big( \|\xi^1\|_{\dbL^{2,\k}_{H}} \!\!+\|\xi^2\|_{\dbL^{2,\k}_{H}} \!\!+ (\phi^{2,\k}_{H})^{1/2}\!\!\Big).
\end{array}
\eea
\end{thm}

\proof (i) First, by Lemma \ref{lem-BSDEest} we have:
\beaa
|\cY^\dbP_t(1,\xi)|^2 \le C_\k\Big(\dbE^\dbP_t\Big[ |\xi|^\k + \int_t^1 |\hat F^0_s|^\k ds\Big]\Big)^{2/\k}, ~\dbP-\mbox{a.s. for all}~\dbP\in{\cP_H^\k}, t\in [0,1].
\eeaa
By the representation formula \reff{representation}, this provides
\beaa
|Y_t|^2 \le C_\k\Big(\dbE^{H,\dbP}_t\Big[ |\xi|^\k + \int_t^1 |\hat F^0_s|^\k ds\Big]\Big)^{2/\k}, ~\dbP-\mbox{a.s. for all}~\dbP\in{\cP_H^\k}, t\in [0,1],
\eeaa
and, by the definition of the norms, we get
\bea
\label{Yest}
\|Y\|^2_{\dbD^{2,\k}_H} \le C_\k\left(\|\xi\|^2_{\dbL^{2,\k}_{H}} + \phi^{2,\k}_{H}\right).
\eea
Next, under each $\dbP\in{\cP_H^\k}$, applying It\^o's formula to $|Y|^2$, it follows from the Lipschitz conditions \reff{FLip-yz} that:
\beaa
 \dbE^\dbP\Big[\int_0^1 |\hat a^{1\slash 2}_s Z_s|^2 ds \Big]
 &\le&
 \dbE^\dbP\Big[|Y_0|^2+ \int_0^1 |\hat a^{1\slash 2}_s Z_s|^2 ds \Big]\\
 &\le&
 C\dbE^\dbP\Big[|\xi|^2 + \int_0^1 |Y_t|\big(|\hat F^0_t| + |Y_t| + |\hat a^{1\slash 2}_t Z_t|\big)ds + \int_0^1 |Y_t| d K^\dbP_t\Big]\\
 &\le&
 C\e^{-1} \dbE^\dbP\Big[|\xi|^2 + \sup_{0\le t\le 1}|Y_t|^2 + \Big(\int_0^1|\hat F^0_t| dt\Big)^{2}\Big]
 \\
 && + \e\dbE^\dbP\Big[\int_0^1|\hat a^{1\slash 2}_t Z_t|^2dt + |K^\dbP_1|^2\Big]
\eeaa
for any $\e\in(0,1]$. By the definition of $K^\dbP$, one gets immediately that
\bea
\label{Kest}
 \dbE^\dbP[|K^\dbP_1|^2]
 &\le& C_0\dbE^\dbP\Big[|\xi|^2 + \sup_{0\le t\le 1}|Y_t|^2
                        + \int_0^1|\hat a_t^{1\slash 2}Z_t|^2dt + \Big(\int_0^1|\hat F^0_t|dt\Big)^{2}
                   \Big],
\eea
for some constant $C_0$ independent of $\e$. Then,
\beaa
 \dbE^\dbP\Big[\int_0^1 |\hat a^{1\slash 2}_s Z_s|^2 ds \Big]
 &\le&
 C\e^{-1} \dbE^\dbP\Big[|\xi|^2 + \sup_{0\le t\le 1}|Y_t|^2 + \Big(\int_0^1|\hat F^0_t|dt\Big)^{2}\Big]
 \\
 &&+ (1+C_0)\e\;\dbE^\dbP\Big[\int_0^1|\hat a^{1\slash 2}_t Z_t|^2dt\Big],
\eeaa
where we recall that the constant $C$ changes from line to line.
By setting $\e := [2(1+C_0)]^{-1}$, this provides
\beaa
\dbE^\dbP\Big[\int_0^1 |\hat a^{1\slash 2}_s Z_s|^2 ds \Big] \le C \dbE^\dbP\Big[|\xi|^2 + \sup_{0\le t\le 1}|Y_t|^2 + \Big(\int_0^1|\hat F^0_t|dt\Big)^2\Big].
\eeaa
By \reff{Yest} and noting that $\phi^{2,1}_{H}\le \phi^{2,\k}_{H}$ for $\k>1$, we have
\bea
\label{Zest}
\|Z\|_{\dbH^{2,\k}_H}^2
\le
C\big(\|\xi\|^2_{\dbL^{2,\k}_{H}} + \phi^{2,\k}_{H}\big).
\eea
This, together with \reff{Yest} and \reff{Kest}, proves \reff{est1}.

\bs

\no (ii) First, following the same arguments as in Lemma \ref{lem-BSDEest}, we have
\beaa
 \big|\cY^\dbP_t(1,\xi_1) - \cY^\dbP_t(1,\xi_2)\big|
 &\le&
 C\Big(\dbE^\dbP_t\Big[|\d \xi|^\k\Big]\Big)^{2/\k},~\dbP-\mbox{a.s. for all}~\dbP\in{\cP_H^\k}, t\in [0,1].
\eeaa
Then, following similar arguments as in (i) we have
\bea
\label{dYest}
\|\d Y\|_{\dbD^{2,\k}_H} \le C\|\d\xi\|_{\dbL^{2,\k}_{H}}.
\eea
Next, under each $\dbP\in{\cP_H^\k}$, applying It\^o's formula to $|\d Y|^2$ we get
\beaa
 \dbE^\dbP\Big[\int_0^1 |\hat a^{1\slash 2}_s \d Z_s|^2 ds \Big]
 &\le&
 \dbE^\dbP\Big[|\d Y_0|^2+\int_0^1 |\hat a^{1\slash 2}_s \d Z_s|^2 ds \Big]
 \\
 &\le& C\dbE^\dbP\Big[|\d \xi|^2
                      + \int_0^1 |\d Y_t|\big(|\d Y_t| + |\hat a^{1\slash 2}_t\d Z_t|\big)ds
                      + \Big|\int_0^1 \d Y_t d (\d K^\dbP_t)\Big|
                 \Big]
 \\
 &\le&
 C\dbE^\dbP\Big[|\d \xi|^2 + \sup_{0\le t\le 1}|\d Y_t|^2 + \sup_{0\le t\le 1}|\d Y_t|[K^{1, \dbP}_1 + K^{2,\dbP}_1]\Big]
 \\
 &&
 + \frac12\dbE^\dbP\Big[\int_0^1|\hat a^{1\slash 2}_t \d Z_t|^2dt\Big].
\eeaa
Then, by \reff{dYest} and \reff{est1},
\beaa
\dbE^\dbP\Big[\int_0^1 |\hat a^{1\slash 2}_s \d Z_s|^2 ds \Big]
&\le& C_\k\|\d\xi\|^2_{\dbL^{2,\k}_{H}} + C_\k\|\d\xi\|_{\dbL^{2,\k}_{H}}\Big(\dbE^\dbP\big[|K^{1, \dbP}_1|^2 + |K^{2,\dbP}_1|^2\big]\Big)^{1/2}\\
&\le& C_\k\|\d\xi\|^2_{\dbL^{2,\k}_{H}} + C_\k\|\d\xi\|_{\dbL^{2,\k}_{H}}\Big(\|\xi^1\|_{\dbL^{2,\k}_{H}}+\|\xi^2\|_{\dbL^{2,\k}_{H}} + (\phi^{2,\k}_{H})^{1/2}\Big).
\eeaa
The estimate for $\d K^\dbP$ is obvious now.
\ep

\quad

We are now ready to state the main result of this paper. Recall that $\cL_{H}^{2,\k}$ is the closure of UC$_b(\O)$ under the norm $\|.\|_{\dbL^{2,\k}_{H}}$.

\begin{thm}
\label{thm-existence}
Let Assumptions \ref{assum-H} and \ref{assum-F} hold. Then for any $\xi\in \cL_{H}^{2,\k}$,  the 2BSDE \reff{2BSDE} has a unique solution $(Y,Z)\in\dbD^{2,\k}_H \times \dbH^{2,\k}_H$.
\end{thm}

\proof (i) We first assume $\xi\in {\rm UC}_b(\O)$.
In this case, by Step 2 of the proof of Theorem 4.5 in \cite{STZ09c}, there exist $(Y, Z)\in \dbD^{2,\k}_H\times \dbH^{2,\k}_H$ such that $Y_1=\xi$, ${\cP_H^\k}-$q.s. and the $K^\dbP$ defined by \reff{KP} is nondecreasing, $\dbP-$a.s. More precisely, $Y_t=V^+_t:=\lim_{\dbQ\ni r\downarrow t}V_r$, where $V$ is defined in that paper.
We notice that the modification of the space of measure $\cP^\k_H$ does not alter the arguments. Moreover, by Proposition 4.10 in \cite{STZ09c}, the representation \reff{representation} holds:
\bea
\label{representation existence}
Y_t = \esup_{\dbP' \in {\cP_H^\k}(t+, \dbP)}~^{\!\!\!\!\!\!\dbP}~ \cY^{\dbP'}_{t}(1,\xi),~~\dbP-\mbox{a.s. for all} ~\dbP\in {\cP_H^\k}, t\in [0,1].
\eea
The construction of $V$ in \cite{STZ09c} is crucially based on the so-called regular conditional probability distribution (r.c.p.d., see Subsection \ref{sect-rcpd}) which allows to define the process $Y$ on $\Omega$ without exception of any zero measure set. Then, $Y$ is shown to satisfy a dynamic programming principle which induces the required decomposition by an appropriate extension of the Doob-Meyer decomposition.

\no (ii) It remains to check the minimum condition \reff{minimum}. We follow the arguments in the proof of Theorem \ref{thm-representation}. For $t\in [0,1]$, $\dbP\in{\cP_H^\k}$, and $\dbP'\in{\cP_H^\k}(t+,\dbP)$, we denote $\d Y:=Y-\cY^{\dbP'}(1,\xi)$, $\d Y:=Z-\cZ^{\dbP'}(1,\xi)$, and we introduce the process $M$ of \reff{M}. Then, it follows from the non-decrease of $K^{\dbP'}$ that
 \bea\label{ineq1existence}
 Y_t-\cY_t^{\dbP'}(1,\xi)
 =
 \d Y_t
 =
 \dbE^{\dbP'}_t\Big[\int_t^1 M_s dK^{\dbP'}_s \Big]
 \ge
 \dbE^{\dbP'}_t\Big[(\inf_{t\le s\le 1} M_t^{-1}M_s) \big(K^{\dbP'}_1-K^{\dbP'}_t\big)\Big].
 \eea

On the other hand, by \reff{Mest} and \reff{ineq1existence}, we estimate by the H\"older inequality that
 \beaa
 &&\dbE^{\dbP'}_t\big[K^{\dbP'}_1-K^{\dbP'}_t\big]\\
 &=&
 \dbE^{\dbP'}_t\Big[\big(\inf_{t\le s\le 1}M_t^{-1} M_s\big)^{1/3}\big(K^{\dbP'}_1-K^{\dbP'}_t\big)^{1/3}
                        \big(\inf_{t\le s\le 1} M_t^{-1}M_s\big)^{-1/3}\big(K^{\dbP'}_1-K^{\dbP'}_t\big)^{2/3}
                   \Big]
 \\
 &\le&
 \Big(\dbE^{\dbP'}_t\big[\big(\inf_{t\le s\le 1}M_t^{-1} M_s\big) \big(K^{\dbP'}_1-K^{\dbP'}_t\big)\big]
\dbE^{\dbP'}_t\big[\sup_{t\le s\le 1} M_tM_s^{-1}\big]\dbE^{\dbP'}_t\big[\big(K^{\dbP'}_1-K^{\dbP'}_t\big)^2\big]\Big)^{1/3}
 \\
 &\le&
 C \Big(\dbE^{\dbP'}_t\big[\big(K^{\dbP'}_1\big)^2\big]
 \dbE^{\dbP'}_t\big[\big(\inf_{t\le s\le 1}M_t^{-1} M_s\big) \big(K^{\dbP'}_1-K^{\dbP'}_t\big)\big]\Big)^{1/3}
 \\
 &\le&
 C \Big(\dbE^{\dbP'}_t\big[\big(K^{\dbP'}_1\big)^2\big]\Big)^{1/3}
 \big(\d Y_t\big)^{1/3}.
 \eeaa
By following the argument of the proof of Theorem \ref{thm-representation} (ii) and (iii), we then deduce that the family $\{K^\dbP,\dbP\in{\cP_H^\k}\}$ inherits the minimum condition \reff{minimum} from \reff{representation existence}.

\ms

(ii) In general, for $\xi\in \cL_{H}^{2,\k}$, by the definition of the space $\cL_{H}^{2,\k}$ there exist $\xi_n\in{\rm UC}_b(\O)$ such that $\lim_{n\to\infty}\|\xi_n - \xi\|_{\dbL^{2,\k}_{H}}=0$. Then it is clear that
\bea
\label{xin}
\sup_{n\ge 1}\|\xi_n\|_{\dbL^{2,\k}_{H}} <\infty &\mbox{and}& \lim_{n, m\to \infty}\|\xi_n - \xi_m\|_{\dbL^{2,\k}_{H}} =0.
\eea
Let $(Y^n, Z^n)\in \dbD^{2,\k}_H\times \dbH^{2,\k}_H$ be the solution to 2BSDE \reff{2BSDE} with terminal condition $\xi_n$, and \bea
\label{KnP}
K^{n,\dbP}_t := Y^n_0 - Y^n_t + \int_0^t \hat F_s(Y^n_s, Z^n_s) ds + \int_0^t Z^n_s dB_s, ~~0\le t\le 1,~~\dbP-\mbox{a.s.}
\eea
By Theorem \ref{thm-apriori}, as $n, m \to\infty$ we have
\beaa
&& \|Y^n-Y^m\|_{\dbD^{2,\k}_H}^2 + \|Z^n-Z^m\|_{\dbH^{2,\k}_H}^2 + \sup_{\dbP\in {\cP_H^\k}}\dbE^\dbP\left[\sup_{0\le t\le 1}|K^{n,\dbP}_t-K^{m, \dbP}_t|^2\right]
\\
&\le&
C_\k\|\xi_n-\xi_m\|_{\dbL^{2,\k}_{H}}^2
+ C_\k\big(\|\xi_n\|_{\dbL^{2,\k}_{H}} + \|\xi_m\|_{\dbL^{2,\k}_{H}}+\|\hat F^0\|_{\dbH^{2,\k}_{H}}\big)
   \|\xi_n-\xi_m\|_{\dbL^{2,\k}_{H}}\to 0.
\eeaa
Then by otherwise choosing a subsequence, we may assume without loss of generality that,
\bea
\label{Ynm}
\|Y^n-Y^{m}\|_{\dbD^{2,\k}_H}^2
+ \|Z^n-Z^{m}\|_{\dbH^{2,\k}_H}^2
+ \sup_{\dbP\in {\cP_H^\k}}\dbE^\dbP\left[\sup_{0\le t\le 1}|K^{n,\dbP}_t-K^{m, \dbP}_t|^2\right]
&\le&
2^{-n},
\eea
for all $m\ge n\ge 1$. This implies that, for every $\dbP\in{\cP_H^\k}$ and $m\ge n\ge 1$,
\bea
\label{Cauchy}
 \dbP\left[\sup_{0\le t\le 1}\big[|Y^n_t-Y^{m}_t|^2
           + |K^{n,\dbP}_t - K^{m,\dbP}_t|^2\big]
           + \int_0^1 |Z^n_t-Z^{m}_t|^2 dt
           > \frac{1}{n}
     \right]
 ~\le~ C n 2^{-n}.
\eea
Define
\bea
\label{limit}
Y := \limsup_{n\to\infty} Y^n, ~~ Z := \limsup_{n\to \infty} Z^n,~~ K^{\dbP} := \limsup_{n\to\infty} K^{n,\dbP},
\eea
where the $\limsup$ for $Z$ is taken componentwise.
It is clear that $Y, Z, K^{\dbP}$ are all $\dbF^+-$progressively measurable. By \reff{Cauchy}, it follows from the Borel-Cantelli Lemma that
\beaa
\lim_{n\to\infty}\Big[\sup_{0\le t\le 1}\big\{|Y^n_t-Y_t|^2 + |K^{n,\dbP}_t - K^{\dbP}_t|^2\big\} + \int_0^1 |Z^n_t-Z_t|^2 dt\Big]= 0, ~~\dbP-\mbox{a.s. for all}~\dbP\in{\cP_H^\k}.
\eeaa
Since $Y^n, K^{n,\dbP}$ are {\cad} and $K^{n,\dbP}$ is nondecreasing, this implies that $Y$ is {\cad}, ${\cP_H^\k}-$q.s. and $K^\dbP$ is {\cad} and nondecreasing, $\dbP-$a.s. Moreover, for every $\dbP\in{\cP_H^\k}$ and $n\ge 1$, sending $m\to\infty$ in \reff{Ynm} and applying Fatou's lemma under $\dbP$,  we obtain:
\beaa
\dbE^{\dbP}\Big[\sup_{0\le t\le 1}\big\{|Y^n_t-Y_t|^2 +|K^{n,\dbP}_t-K^{\dbP}_t|^2\big\} + \int_0^1 |Z^n_t-Z_t|^2\Big] \le 2^{-n}.
\eeaa
This implies that
\beaa
 \|Y^n-Y\|_{\dbD^{2,\k}_H}^2 + \|Z^n-Z\|_{\dbH^{2,\k}_H}^2 + \sup_{\dbP\in {\cP_H^\k}}\dbE^\dbP\left[\sup_{0\le t\le 1}|K^{n,\dbP}_t-K^{\dbP}_t|^2\right]
 &\le&
 2^{-n}\to 0, ~~\mbox{as}~n\to\infty.
\eeaa
Then it is clear that $(Y,Z)\in\dbD^{2,\k}_H \times \dbH^{2,\k}_H$.

Finally, since $(Y^n, Z^n, K^{n,\dbP})$ satisfy \reff{2BSDEP} and  \reff{representation}, the limit $(Y,Z, K^\dbP)$ also satisfies \reff{2BSDEP} and \reff{representation}. Then by the proof of Theorem \ref{thm-existence}, the family $\{K^\dbP,\dbP\in\cP^\k_H\}$ satisfies \reff{minimum}. Hence $(Y,Z)$ is a solution to 2BSDE \reff{2BSDE}.
\ep

\begin{rem}
\label{rem-Nutz}{\rm
After the completion of this paper, Marcel Nutz pointed out that our solution of the 2BSDE in the present contexts is in fact $\dbF-$progressively measurable, as a consequence of the uniform continuity in $\o$ in our setting. See Proposition 4.11 in \cite{STZ09c}. However, the $\dbF^+-$progressive measurability developed in this paper seems to be more robust to potential extensions of the spaces.
}\end{rem}

\section{Connection with fully nonlinear PDEs}
\label{sect-PDE}
\setcounter{equation}{0}

\subsection{The Markovian setup}

In this section we consider the case:
 \beaa
 H_t(\omega,y,z,\g)
 &=&
 h(t,B_t(\omega),y,z,\g),
 \eeaa
where $h:  [0,1]\times \dbR^d \times \dbR\times \dbR^d\times D_h \to \dbR$ is a deterministic map.
Then the corresponding conjugate and bi-conjugate functions become
 \bea
 f(t,x,y,z,a)
 &:=&
 \sup_{\gamma\in D_h}~\big\{\frac{1}{2}a:\g - h(t,x,y,z,\g)\big\}, ~~a\in \dbS^{>0}_d,
 \label{F-Markov}
 \\
 \hat h(t,x,y,z,\g)
 &:=&
 \sup_{a\in \dbS^{>0}_d}~\big\{\frac{1}{2}a:\g - f(t,x,y,z,a)\big\}, ~~\g\in \dbR^{d\times d}.
 \label{hatH}
 \eea
Notice that $-\infty<\hat h\leq h$ and $\hat h$ is nondecreasing convex in $\g$. Also, $\hat h=h$ if and only if $h$ is convex and nondecreasing in $\g$.

In the present context, we write $\cP^\k_h:=\cP^\k_H$.
The following is a slight strengthening of  Assumption \ref{assum-H} to our Markov framework.

\begin{assum}\label{assum-f}
$\cP^\k_h$ is not empty, the domain $D_{f_t}$ of the map $a\longmapsto f(t,x,y,a)$ is independent of $(x,y,z)$. Moreover, on $D_{f_t}$, $f$ is uniformly continuous in $t$, uniformly in $a$, and for some constant $C$ and modulus of continuity $\rho$ with polynomial growth:
\begin{equation}
\label{fLip}
\big|f(t, x,y,z, a) - f(t,x',y',z', a)\big|
\le
\rho(|x-x'|)+C\left(|y-y'|+\big|a^{1\slash 2}(z_1-z_2)\big|\right),
\end{equation}
for all $t\in [0,1]$, $a\in D_{f_t}$, $x,x',z,z'\in\dbR^d$, $y,y'\in\dbR$.
\end{assum}

Next, let  $g: \dbR^d\to \dbR$ be a Lebesgue measurable function. In this section we shall always consider the 2BSDE \reff{2BSDE} in this Markovian setting with terminal condition $\xi = g(B_1)$:
 \bea\label{2BSDEMarkov}
 Y_t = g(B_1) - \int_t^1 f(s,B_s,Y_s, Z_s, \hat a_s) ds - \int_t^1 Z_s dB_s + K_1 - K_t, ~~0\le t\le 1, ~~\mbox{${\cP_H^\k}-$q.s.}
 \eea
Our main objective is to establish the connection $Y_t = v(t,B_t)$, $t\in[0,1]$, ${\cP_H^\k}-$q.s. where $v$ solves, in some sense, the following fully nonlinear PDE:
 \bea
 \label{PDE}
 \left\{\ba{l}
 \cL v (t,x) := \pa_tv(t,x) + \hat h\big(t, x, v(t,x), Dv(t,x), D^2 v(t,x)\big)=0, ~~0\le t<1,
 \\
 v(1,x) = g(x).
 \ea\right.
\eea
We remark that the nonlinearity of the above PDE is the nondecreasing and convex envelope $\hat h$, not the original $h$. This is illustrated by the following example.

\begin{eg}
\label{eg-GammaConstraint1}
{\rm The problem of hedging under gamma constraints in dimension $d=1$, as formulated by Cheridito, Soner and Touzi \cite{CST}, leads to the specification
 \beaa
 h(t,x,y,z,\g):=\frac12\g
 ~\mbox{if}~
 \g\in[\underline{\G},\overline{\G}],
 &\mbox{and}&
 \infty~\mbox{otherwise,}
 \eeaa
where $\underline{\G}<0<\overline{\G}$ are given constants. Then, direct calculation leads to
 \beaa
 f(a)
 &=&
 \frac12 \big(\overline{\G}(a-1)^+-\underline{\G}(a-1)^-\big),
 ~~
 a> 0,
 \\
 \hat{h}(\g)
 &=&
 \frac12(\g\vee\underline{\G})
 ~~\mbox{if}~~
 \g\le\overline{\G},
 ~~\mbox{and}~~
 \infty~~\mbox{otherwise.}
 \eeaa
\ep
}
\end{eg}

We will discuss further this case in Example \ref{eg-GammaConstraint2} below, in order to obtain the nonlinearity appearing in the PDE characterization of \cite{CST} for the superhedging problem under gamma constraints. Indeed, equation \reff{PDE} needs to be reformulated in some appropriate sense if $D_{h}\neq\dbS_d$, because then $\hat h$ may take infinite values, and the meaning of \reff{PDE} is not clear anymore. This leads typically to a boundary layer and requires the interpretation of the equation in the relaxed boundary value sense of viscosity solutions, see, e.g. \cite{cil}.

\subsection{A nonlinear Feynman-Kac representation formula}
\label{sect-PDE2}

\begin{thm}\label{thm-FK}
Let Assumption \ref{assum-f} hold true. Suppose further that $\hat h$ is continuous in its domain, $D_f$ is independent of $t$ and is bounded both from above and away from $0$. Let $v\in C^{1,2}([0,1),\dbR^d)$ be a classical solution of \reff{PDE} with  $\{(v,Dv)(t,B_t),t\in[0,1]\}\in\dbD^{2,\k}_H\times\dbH^{2,\k}_H$. Then:
 \beaa
 &Y_t:=v(t,B_t),~Z_t:=Dv(t,B_t),~K_t:=\int_0^t k_s ds&
 \\
 &\mbox{with}~k_t:=\hat h\left(t,B_t,Y_t,Z_t,\Gamma_t\right)-\frac12\hat a_t \: \Gamma_t+f\left(t,B_t,Y_t,Z_t,\hat a_t\right)
 ~\mbox{and}~\Gamma_t:=D^2v(t,B_t)&
 \eeaa
is the unique solution of the 2BSDE \reff{2BSDEMarkov}.
\end{thm}

\proof
By definition $Y_1=g(B_1)$ and \reff{2BSDEMarkov} is verified by immediate application of It\^o's formula. It remains to prove the minimum condition:
 \bea\label{minimum-FK}
 \einf_{\dbP'\in\cP_H^\k(t+,\dbP)}\dbE_t^{\dbP'}\left[\int_t^1 k_s ds\right]
 =
 0
 &\mbox{for all}&
 t\in[0,1],~\dbP\in\cP_H^\k,
 \eea
by which we can conclude that $(Y,Z,K)$ is a solution of the 2BSDE \reff{2BSDEMarkov}. Since $g(B_1)\in\dbL_H^{2,\k}$, the uniqueness follows from Theorems \ref{thm-representation} and \ref{thm-apriori} (i).

To prove \reff{minimum-FK}, we follow the same argument as in the proof of Lemma 3.1 in \cite{EPQ}.  For every $\e>0$, notice that the set
 \beaa
 A^\eps
 &:=&
 \left\{ a\in D_{f}:~
        \hat h(t,B_t,Y_t,Z_t,\Gamma_t)
        \le
        \frac12 a:\Gamma_t-f(t,B_t,Y_t,Z_t,a)+\e
 \right\}
 \eeaa
is not empty. Then it follows from a measurable selection argument that there exists a predictable process $a^\eps$ taking values in $D_{f}$ such that
 \beaa
\hat h(t,B_t,Y_t, Z_t,\Gamma_t) &\le& \frac12 a^\e_t:\Gamma_t-f(t,B_t,Y_t,Z_t,a^\e_t)+\e.
 \eeaa
We note that this in particular implies that $\G_t\in D_{\hat h}$.

In the remainder of this proof, we show the existence of an $\dbF-$progressively measurable process $\alpha^\eps$ with values in $\dbS_d^{>0}$ and $\int_0^1|\alpha^\eps_s|ds<\infty$ such that, $\dbP^{\alpha^\eps}-$a.s., $\hat{a}$ is in $A^\eps$. We recall from Remarks \ref{rem-a1} and \ref{rem-a2} that this is not guaranteed in general. Notice that this technical difficulty is inherent to the problem and requires to be addressed even if a maximizer for $\hat h$ does exist.

Let $\dbP:=\dbP^\a\in \cP_H$  and $t_0\in [0,1]$ be fixed. Let
 \beaa
 \t^\e_0
 &:=&
 1\wedge\inf\left\{t\ge t_0\ |\ K_t\ge K_{t_0}+\eps \right\},
 \eeaa
and define:
 \beaa
 \t^\e_{n+1} :=
 1\wedge\inf\Big\{t\ge \t^\e_n\ |\ \hat h(t,B_t,Y_t,\Gamma_t)
                                &\ge&
                                \frac12 a^\e_{\t^\e_n}:\Gamma_t
                                -f(t,B_t,Y_t,Z_t,a^\e_{\t^\e_n})
                                +2\e
            \Big\},
 \eeaa
for $n\ge 0$. Since $K$ is continuous, notice that $\tau^\eps_0>t_0$,~$\cP^\k_H-$q.s.. Also, since $B, Y, Z, \G$ are all continuous in $t$, $\tau^\e_n$ are $\dbF-$stopping times and, for any fixed $\o$, are uniformly continuous in $t$.

Next, for any fixed $a\in D_f$, the function $f(.,a)$ is continuous. Also $\hat h$ is continuous. Then for $\cP_H^\k-$q.s. $\o\in \O$,
 \beaa
 \hat h(t,B_t(\o),Y_t(\o), Z_t(\o),\Gamma_t(\o)) - \frac12 a^\e_{\t^\e_n}(\o) :\Gamma_t(\o)+f(t,B_t(\o),Y_t(\o),Z_t(\o),a^\e_{\t^\e_n}(\o))
 \eeaa
 is uniformly continuous in $t$ for $t\in [\t^\e_n(\o), 1]$. Then $\t^\e_{n+1}(\o) - \t^\e_n(\o)\ge \d(\e,\o)>0$ whenever $\t^\e_{n+1}(\o)<1$, where the constant $\d(\e,\o)$ does not depend on $n$. This implies that $\t^\e_n(\o) = 1$ for $n$ large enough. Applying the arguments in Example 4.5 of \cite{STZ09a} on $[\tau^\eps_0,1]$, one can easily see that there exists an $\dbF-$progressively measurable process $\a^\e$ taking values in $D_f$ such that
 \beaa
 \a^\e_t = \a_t ~~\mbox{for}~ t\in [0,\tau^\eps_0]
 &\mbox{and}&
 \hat a_t = \sum_{n=0}^\infty a^\e_{\t^\e_n} \1_{[\t^\e_n, \t^\e_{n+1})}(t), dt\times d\dbP^{\a^\e}-\mbox{a.s. on}
 ~[\tau^\eps_0,1]\times \O.
 \eeaa
 This implies that
 \beaa
 \hat h(t,B_t,Y_t, Z_t,\Gamma_t) \le \frac12 \hat a_t:\Gamma_t-f(t,B_t,Y_t,Z_t,\hat a_t)+2\e,~dt\times d\dbP^{\a^\e}-\mbox{a.s. on}~[\tau^\eps_0,1]\times \O,
 \eeaa
Under our conditions it is obvious that $\dbP^{\a^\e}\in \cP^\k_H$, then $\dbP^{\a^\e}\in \cP^\k_H(t_0+,\dbP)$ sinse $\tau^\eps_0>t_0$.
Therefore,
 \beaa
 \einf_{\dbP'\in\cP^\k_H(t_0+,\dbP)}~^{\!\!\!\!\!\!\dbP}~\dbE_{t_0}^{\dbP'}\left[\int_{t_0}^1 k_t dt\right]
 &\le&
 \e+\dbE_{t_0}^{\dbP^{\a^\e}}\left[\int_{\tau^\eps_0}^1 k_t dt\right]
 \;\le\;
 \e+2\eps(1-t_0),~~\dbP-\mbox{a.s.}
 \eeaa
By the arbitrariness of $\e>0$, and the nonnegativity of $k$, this provides \reff{minimum-FK}.
\ep

\subsection{Markovian solution of the 2BSDE}

Following the classical terminology in the BSDE literature, we say that the solution of the 2BSDE is Markovian if it can be represented by means of a determinitic function of $(t,B_t)$.
In this subsection we construct a deterministic function $u$, by using a probabilistic representation in the spirit of \reff{representation}, and show its connection with 2BSDE \reff{2BSDEMarkov}. The connection between $u$ and the PDE \reff{PDE} will be established in the next subsection.

Following \cite{STZ09c}, we introduce the shifted probability spaces.
For $0\le t\le 1$, denote by $\O^{t}:= \{\o\in C([t,1], \dbR^d): \o(t)=0\}$ the shifted canonical space; $B^{t}$ the shifted canonical process on $\O^{t}$; $\dbP^{t}_0$ the shifted Wiener measure; $\dbF^{t}$ the shifted filtration generated by $B^{t}$,
$\overline\cP^t_S$ the corresponding collection of martingale measures induced by the strong formulation, and $\hat a^t$ the universal quadratic variation density of $B^t$.  In light of Definition \ref{defn-cP}, we define

\begin{defn}
\label{defn-cPt}
For $t\in [0,1]$, let $\cP^{\k,t}_h$ denote the collection of all those $\dbP\in \overline\cP^t_S$ such that
 \bea
 \label{ellipticity-t}
 \ba{c}
 \dis \underline a_\dbP \le \hat a^t \le \overline a_\dbP, ~~
 ds\times d\dbP-\mbox{a.s. on}~[t,1]\times\O^t,~~
 \mbox{for some}~ \underline a_\dbP, \overline a_\dbP\in \dbS^{>0}_d,
 \\
 \dis \mbox{and}~~
 \dbE^{\dbP}\Big[\big(\int_t^1 |\hat f^{t,0}_s)|^\k ds\big)^{2/\k}
            \Big]
 <
 \infty ,~~\mbox{where}~~\hat f^{t,0}_s:= f(s, 0, 0,0,\hat a^t_s).
\ea\eea
\end{defn}

\begin{rem}
{\rm By Lemma \ref{lem-rcpd} below, $\cP^\k_h\neq\emptyset$ implies that $\cP^{\k,t}_h\neq\emptyset$ for all $t\in[0,1]$.
\ep}
\end{rem}

By Assumption \ref{assum-f}, the polynomial growth of $\rho$, and the first part of \reff{ellipticity-t}, it is clear that
 \beaa
 \dbE^{\dbP}\Big[\big(\int_t^1 |\hat f^{t,0}_s|^\k ds\big)^{2/\k} \Big]
 < \infty ~
 \mbox{if and only if}~
 \dbE^{\dbP}\Big[\big(\int_t^1 |f(s, B^t_s, 0,0,\hat a^t_s)|^\k ds
                 \big)^{2/\k}
            \Big]
 <\infty,
\eeaa
and thus, for $t=0$, we see that $\cP^\k_h = \cP^{\k,0}_h$ as defined in Definition \ref{defn-cP}.

We next define a similar notation to \reff{BSDEP}. For any $(t,x)\in[0,1]\times\dbR^d$,  denote
 \beaa
 B^{t,x}_s := x+B^t_s  &\mbox{for all}& s\in[t,1].
 \eeaa
Let $\t$ be and $\dbF^t-$stopping time, $\dbP\in \cP^{\k,t}_h$, and $\eta$ a $\dbP-$square inetgarble $\cF^t_\t-$measurable r.v. See Remark \ref{rem-FF+}. We denote by $\big(\cY^{\dbP}, \cZ^{\dbP}\big):=\big(\cY^{t,x,\dbP}(\t,\eta), \cZ^{t,x,\dbP}(\t, \eta)\big)$ the solution of the following BSDE:
\bea
\label{BSDEt}
\cY^\dbP_s = \eta - \int_t^{\t} f(r, B^{t,x}_r, \cY^{\dbP}_r, \cZ^{\dbP}_r, \hat a_r) dr
                  - \int_s^{\t} \cZ^{\dbP}_r dB_r,~t\le s\le \t,~\dbP-\mbox{a.s.}
\eea
Similar to \reff{BSDEP}, under our assumptions the above BSDE has a unique solution. We now introduce the value function:
\bea
\label{u}
u(t,x) := \sup_{\dbP \in \cP^{\k,t}_h} \cY^{t,x,\dbP}_t\big(1, g(B_1^{t,x})\big),
&\mbox{for}&
(t,x)\in[0,1]\times\dbR^d.
\eea
By the Blumenthal zero-one law \reff{0-1 MRP}, it follows that $\cY^{t,x,\dbP}_t\big(1, g(B^{t,x}_1)\big)$ is a constant and thus $u(t,x)$ is deterministic.

\begin{rem}{\rm
Notice that, in contrast with the previous sections, we are now implicitly working with the filtration $\dbF$. However, the subsequent Theorem \ref{thm-2BSDEu} connects $u(t,B_t)$ to the solution of the 2BSDE, implying that $Y$ is  $\dbF-$progressively measurable. See Remark \ref{rem-Nutz}.
}
\end{rem}

We next state a strengthening of Assumption \ref{assum-F} in the present Markov framework.
\begin{assum}\label{assum-g}
The function $g$ has polynomial growth, and there exists a continuous positive function $\L(t,x)$ such that, for any $(t,x)$:
 \bea
 \label{Ldominate}
 &&\sup_{\dbP\in\cP^{\k,t}_h}\dbE^{\dbP}\Big[|g(B^{t,x}_{1})|^\k + \int_t^{1} |f(s, B^{t,x}_s, 0,0, \hat a^t_s)|^\k ds\Big]\le \L^\k(t,x),
 \\
 \label{Lintegrability}
  &&\sup_{\dbP\in\cP^{\k,t}_h}\dbE^{\dbP}\Big[\sup_{t\le s\le 1}\L^2(s,B^{t,x}_s)\Big] <\infty.
 \eea
\end{assum}
By the definition of $\L$, it is clear that 
\bea
\label{u<Lamda}
|u| &\le& \L.
\eea

\begin{rem}
\label{rem-g}
{\rm There are two typical sufficient conditions for the existence of such $\L$:
\\
(i) $f$ and $g$ are bounded. In this case one can choose $\L$ to be a constant.
\\
(ii) $D_f$ is bounded and $\sup_{\dbP\in\cP^{\k,t}_h}
\dbE^{\dbP}\Big[\int_t^{1} |\hat f^{t,0}_s|^\k ds\Big] \le C$ for all $t$. In this case one can choose $\L$ to be a polynomial of $|x|$.
\ep}
\end{rem}

\begin{thm}
\label{thm-2BSDEu}
Let Assumptions \ref{assum-f} and \ref{assum-g} hold true, and $g$ be uniformly continuous, so that the 2BSDE \reff{2BSDEMarkov} has a unique solution $(Y, Z)\in \dbD^{2,\k}_H\times \dbH^{2,\k}_H$. Then $Y_t = u(t,B_t)$. Moreover, $u$ is uniformly continuous in $x$, uniformly in $t$, and right continuous in $t$.
\end{thm}

\proof The wellposedness of 2BSDE \reff{2BSDEMarkov} follows directly from Theorem \ref{thm-existence}. Notice that $u(t, B_t)=V_t$ as defined in \cite{STZ09c}. By Remark \ref{rem-Nutz}, $Y_t = V_t$, and thus $Y_t = u(t, B_t)$.

The uniform continuity of $u$ follows from Lemma 4.6 of \cite{STZ09c}; alternatively one can follow the proof of Lemma \ref{lem-BSDEest} applied to the difference of two solutions. Finally, for any $(t,x)$ and $\d>0$, the decomposition 
\beaa
|u(t+\d, x)- u(t,x)| = u(t+\d, x)- u(t+\d, B^{t,x}_{t+\d}) +  Y^{t,x}_{t+\d} - Y^{t,x}_t
\eeaa
implies the right continuity of $u$ in $t$, as a consequence of the uniform continuity of $u$ in $x$, uniformly in $t$,  and the right continuity of the process $Y$.
\ep

\vspace{5mm}

Finally, for later use, we provide an additional regularity result on $u$.

\begin{prop}
\label{prop-u}
Let Assumptions \ref{assum-f} and \ref{assum-g}  hold true, and  $g$ be lower-semicontinuous. Then $u$ is lower-semicontinuous in $(t,x)$.
\end{prop}
The proof is closely related to the Dynamic Programming Principle, and is postponed to Subsection \ref{subsect-DPlsc}.

\subsection{The viscosity solution property}

We shall make use of the classical notations in the theory of viscosity solutions:
 \bea
 \label{h*}
 \ba{lll}
 \dis u_*(\th):=\liminf_{\th'\to \th}u(\th)
 &\mbox{and}&
 u^*(\th):=\limsup_{\th'\to \th}u(\th'), ~~\mbox{for}~~\th = (t,x);
 \\
 \dis \hat h_*(\theta) :=\liminf_{\theta'\to \theta} \hat h(\theta')
 &\mbox{and}&
 \hat h^*(\theta) :=\limsup_{\theta'\to\theta} \hat h(\theta'),~~\mbox{for}~~\th = (t,x,y,z,\g).
 \ea
 \eea

\begin{thm}
\label{thm-viscosity}
Let Assumptions \ref{assum-f} and \ref{assum-g} hold true. Then:
\\
{\rm (i)}\q $u$ is a viscosity subsolution of
\bea
\label{PDE2}
-\pa_t u^* - \hat h^*(\cd,u^*, Du^*, D^2u^*) \le 0
&\mbox{on}&
[0,1)\times\dbR^d.
\eea
{\rm (ii)}\q Assume further that $g$ is lower-semicontinuous and $D_f$ is independent of $t$, then $u$ is a viscosity supersolution of
\bea
\label{PDE1}
-\pa_t u_* - \hat h_*(\cd,u_*, Du_*, D^2u_*) \ge 0
&\mbox{on}&
[0,1)\times\dbR^d.
\eea
\end{thm}

\begin{eg}
\label{eg-GammaConstraint2}
{\rm Let us illustrate the role of $\hat h^*$ and $\hat h_*$ in the context of Example \ref{eg-GammaConstraint1}. In this case, one can check immediately that
\beaa
\hat{h}_*
 = \hat h
 &\mbox{and}&
 \hat{h}^*(\g)
 =
 \frac12(\g\vee\underline{\G})\1_{\{\g < \overline{\G}\}} + \infty\1_{\{\g \ge \overline{\G}\}}.
 \eeaa
 Then the above viscosity properties are equivalent to
 \beaa
 \min\big\{ -\pa_t u^* - \frac12(D^2u^* \vee\underline{\G}), ~ \bar\G - D^2u^*\big\} \le 0,
 \\
 \min\big\{ -\pa_t u_* - \frac12(D^2u_* \vee\underline{\G}), ~ \bar\G - D^2u_*\big\} \ge 0,
 \eeaa
which is exactly the nonlinearity obtained in \cite{CST}.
 \ep
}
\end{eg}

\begin{rem}
\label{rem-viscosity}
{\rm
(i) If $u$ is continuous and $D_{\hat h}=\dbR^{d\times d}$, then by Theorem \ref{thm-viscosity} $u$ is a viscosity solution to PDE \reff{PDE} in the standard sense.
\\
\noindent (ii) If the comparison principle for the following relaxed boundary value fully nonlinear PDE \reff{PDE2}-\reff{PDE1} with boundary condition holds:
 \begin{equation}\label{relaxedboundary}
 \begin{array}{c}
 \max\left\{ \big(-\pa_t v - \hat h_*(\cd,v, Dv, D^2v)\big)(T,.), v(T,.)-g\right\} \;\ge\; 0
 \\
 \min\left\{ \big(-\pa_t v - \hat h^*(\cd,v, Dv, D^2v)\big)(T,.), v(T,.)-g\right\} \;\le\; 0
 \end{array}
 \end{equation}
then $u$ is continuous and is the unique viscosity solution to the above problem.  We refer to Crandal, Ishii and Lions \cite{cil} for the notion of relaxed boundary problems.
 \ep}
\end{rem}

The viscosity property is a consequence of the following dynamic programming principle.

\begin{prop}\label{prop-DP}
Let $g$ be lower-semicontinuous, $t\in [0,1]$, and $\{\t^\dbP, \dbP\in \cP^{\k,t}_h\}$ be a family of $\dbF^t-$stopping times. Then, under Assumptions \ref{assum-f} and \ref{assum-g}:
 \beaa
 u(t,x)
 &=&
 \sup_{\dbP \in \cP^{\k,t}_h} \cY^{t,x,\dbP}_t\big(\tau^\dbP, u(\tau^\dbP,B_{\tau^\dbP}^{t,x})\big).
 \eeaa
\end{prop}

The proof of Proposition \ref{prop-DP} is reported in subsections \ref{sect-weakDP1} and \ref{subsect-DPlsc}.

\bs
\no{\it Proof of Theorem \ref{thm-viscosity}.} (i) We argue by contradiction, and we aim for a contradiction of the dynamic programming principle. Assume to the contrary that
 \bea\label{maxstrict}
 0=(u^*-\varphi)(t_0,x_0)>(u^*-\varphi)(t,x)
 &\mbox{for all}&
 (t,x)\in([0,1]\times\dbR^d)\setminus\{(t_0,x_0)\}
 \eea
for some $(t_0,x_0)\in[0,1)\times\dbR^d$ and
 \bea
 \label{Lf>0}
 \big(-\partial_t\varphi-\hat h^*(.,\varphi,D\varphi,D^2\varphi)\big)(t_0,x_0)>0,
 \eea
for some smooth function $\f$. By \reff{u<Lamda}, without loss of generality we may assume $|\f|\le \L$. We note that \reff{Lf>0} implies that $D^2\f(t_0,x_0)\in D_{\hat h}$. Since $\hat h^*$ is upper-semicontinuous and $\f$ is smooth, there exists an open ball $O_r(t_0,x_0)$,  centered at $(t_0,x_0)$ with radius $r$, such that
 \beaa
 -\partial_t\varphi-\hat h(.,\varphi,D\varphi,D^2\varphi)\ge 0, &\mbox{on}& O_r(t_0,x_0).
 \eeaa
Then, we deduce from the definition of $\hat h$  that
 \bea\label{pde>0}
 -\partial_t\varphi-\frac12\alpha\:D^2\varphi+f(.,\varphi,D\varphi,\alpha)\ge 0
 ~\mbox{on}~O_r(t_0,x_0)
 &\mbox{for all}&
 \alpha\in\dbS_d^{>0}(\dbR).
 \eea

By the strict maximum property \reff{maxstrict}, we notice that
 \bea\label{etasupersol}
 \eta \;:=\; -\max_{\partial O_r(t_0,x_0)} (u^*-\varphi) &>& 0.
 \eea
 Let $(t_n,x_n)$ be a sequence of $O_r(t_0,x_0)$ such that
 \beaa
 (t_n,x_n)\longrightarrow (t_0,x_0)
 &\mbox{and}&
 u(t_n,x_n)\longrightarrow u^*(t_0,x_0),
 \eeaa
and define the stopping time $\tau_n:=\inf\{ s>t_n: (s,B^{t_n,x_n}_s)\not\in O_r(t_0,x_0)\}$. Without loss of generality we may assume $r< 1-t_0$, then $\t_n< 1$ and thus $(\t_n,B^{t_n,x_n}_{\t_n})\in \partial O_r(t_0,x_0)$. With this construction we have
 \bea\label{cntaun1}
 c_n:=(\varphi-u)(t_n,x_n)\to 0
 &\mbox{and}&
 u^*(\tau_n,B_{\tau_n}^{t_n,x_n})\le \varphi(\tau_n,B_{\tau_n}^{t_n,x_n})-\eta,
 \eea
by the continuity of the coordinate process.

For any  $\dbP^n\in \cP^{\k,t_n}_h$, we now compute by the comparison result for BSDEs and classical estimates that
 \bea\label{contraDP1}
 \begin{array}{rcl}
 \cY_{t_n}^{t_n, x_n,\dbP^n}\big(\tau_n, u^*(\tau_n, B^{t_n,x_n}_{\tau_n})\big) &-& u(t_n,x_n)
 \\
 &\le& \cY_{t_n}^{t_n, x_n,\dbP^n}\big(\tau_n, \varphi(\tau_n, B^{t_n,x_n}_{\tau_n})-\eta\big)
  -\varphi(t_n,x_n)+c_n
 \\
 &\le&
 \cY_{t_n}^{t_n, x_n,\dbP^n}\big(\tau_n, \varphi(\tau_n, B^{t_n,x_n}_{\tau_n})\big)
  -\varphi(t_n,x_n)+c_n - \eta'
 \end{array}
 \eea
for some positive constant $\eta'$ independent of $n$. Set
 \beaa
 (Y^n,Z^n):=\big(\cY^{t_n,x_n,\dbP^n},\cZ^{t_n,x_n,\dbP^n}\big)\big(\tau_n, \varphi(\tau_n, B^{t_n,x_n}_{\tau_n})\big),
 \eeaa
 \beaa
 \delta Y^n_s:= Y^n_s-\varphi(s,B^{t_n,x_n}_s),
 ~\mbox{and}~
 \delta Z^n_s:= Z^n_s-D\varphi(s,B^{t_n,x_n}_s).
 \eeaa
It follows from It\^o's formula together with the Lipschitz properties of $f$ that, $\dbP^n-$a.s.
 \beaa
 d(\delta Y^n_{s})
 &=&
 \big(-\partial_t\varphi-\frac12\hat a_s\:D^2\varphi+f(.,Y^n_s,Z^n_s,\hat a_s)\big)(s,B^{t_n,x_n}_s)ds+\delta Z^n_sdB_s
 \\
 &=&
 (\phi^n_s+\lambda_s\delta Y^n_s+\delta Z^n_s\bar\alpha^{1/2}\beta_s)ds +\delta Z^n_sdB_s
 \eeaa
where $\lambda$ and $\beta$ are bounded progressively measurable processes, and
 \beaa
 \phi^n_s := \big(-\partial_t\varphi-\frac12\hat a_s\:D^2\varphi+f(.,\varphi,D\varphi,\hat a_s)\big)(s,B^{t_n,x_n}_s)
 \ge
 0 &\mbox{for}& s\in [t_n, \t_n],
 \eeaa
by \reff{pde<0} and the definition of $\t_n$. Let $M$ be defined by \reff{M}, but starting from $t_n$ and under $\dbP^n$.
Then
 \beaa
 \cY_{t_n}^{t_n, x_n,\dbP^n}\big(\tau_n, \varphi(\tau_n, B^{t_n,x_n}_{\tau_n})\big)-\varphi(t_n,x_n)
 =\delta Y^n_{t_n}
 \le
 \dbE^{\dbP^n}\left[M_{\t_n}\delta Y^n_{\t_n}\right]=0.
 \eeaa
Plugging this in \reff{contraDP1}, we get
 \beaa
 \cY_{t_n}^{t_n, x_n,\dbP^n}\big(\tau_n, u^*(\tau_n, B^{t_n,x_n}_{\tau_n})\big) - u(t_n,x_n)
 &\le&
 c_n-\eta'.
 \eeaa
 Note that $\dbP^n\in \cP^{\k, t_n}_h$ is arbitrary and $c_n$ does not depend on $\dbP^n$. Then
 \beaa
 \sup_{\dbP\in \cP^{\k,t_n}_h} \cY_{t_n}^{t_n, x_n,\dbP}\big(\tau_n, u^*(\tau_n, B^{t_n,x_n}_{\tau_n})\big) - u(t_n,x_n)
 &\le&
 c_n-\eta' < 0,
 \eeaa
for large $n$.
This is in contradiction with the dynamic programming principle of Proposition \ref{prop-DP} (or, more precisely, Lemma \ref{lem-bt} below to avoid the condition that $g$ is lower-semicontinuous) .

\ms

\noindent (ii) \quad We again argue by contradiction, aiming for a contradiction of the dynamic programming principle of Proposition \ref{prop-DP}. Assume to the contrary that
 \bea\label{minstrict}
 0=(u_*-\varphi)(t_0,x_0)<(u_*-\varphi)(t,x)
 &\mbox{for all}&
 (t,x)\in([0,1]\times\dbR^d)\setminus\{(t_0,x_0)\}
 \eea
for some $(t_0,x_0)\in[0,1)\times\dbR^d$ and
 \beaa
 \big(-\partial_t\varphi-\hat h_*(.,\varphi,D\varphi,D^2\varphi)\big)(t_0,x_0)<0,
 \eeaa
 for some smooth function $\f$. By \reff{u<Lamda}, without loss of generality we may assume again that $|\f|\le \L$. Note that $\hat h_* \le \hat h$. Then
 \beaa
 \big(-\partial_t\varphi-\hat h(.,\varphi,D\varphi,D^2\varphi)\big)(t_0,x_0)<0.
 \eeaa
 If $D^2\f (t_0, x_0) \in D_{\hat h}$, then it follows from the definition of $\hat h$ that
 \bea
 \label{Lf<0}
 \Big(-\partial_t\varphi-\frac12\bar\alpha\:D^2\varphi+f(.,\varphi,D\varphi,\bar\alpha)\Big)(t_0,x_0) < 0
 \eea
 for some $\bar\a \in \dbS^{>0}_d$. In particular, this implies that $\bar\a\in D_f$. If $D^2\f (t_0, x_0) \notin D_{\hat h}$, since $\pa_t\f(t_0,x_0)$ is finite, we still have $\bar\a \in D_f$ so that \reff{Lf<0} holds. Now by the smoothness of $\f$ and \reff{fLip}, and recalling that $D_f$ is independent of $t$, there exists an open ball $O_r(t_0,x_0)$ with $0<r<1-t_0$ such that
 \bea
 \label{pde<0}
 -\partial_t\varphi-\frac12\bar\alpha\:D^2\varphi+f(.,\varphi,D\varphi,\bar\alpha)\le 0
 ~\mbox{on}~O_r(t_0,x_0).
 \eea

By the strict minimum property \reff{minstrict}, we notice that
 \bea\label{etasubsol}
 \eta \;:=\; \min_{\partial B_r(t_0,x_0)} (u_*-\varphi) &>& 0.
 \eea
As in (i), we consider a sequence $(t_n,x_n)$ of $O_r(t_0,x_0)$ such that
 \beaa
 (t_n,x_n)\longrightarrow (t_0,x_0)
 &\mbox{and}&
 u(t_n,x_n)\longrightarrow u_*(t_0,x_0),
 \eeaa
and we define the stopping time $\tau_n:=\inf\{ s>t_n: (s,B^{t_n,x_n}_s)\not\in O_r(t_0,x_0)\}$, so that
 \bea\label{cntaun2}
 c_n:=(u-\varphi)(t_n,x_n)\to 0
 &\mbox{and}&
 u_*(\tau_n,B_{\tau_n}^{t_n,x_n})\ge \varphi(\tau_n,B_{\tau_n}^{t_n,x_n}) + \eta.
 \eea
For each $n$, let $\bar\dbP^n:=\dbP^{\bar\alpha}\in \overline\cP_S^{t_n}$ be the local martingale measure induced by the constant diffusion $\bar\alpha$. By \reff{fLip}, one can easily see that $\bar \dbP^n\in {\cP_H^{\k,t_n}}$. We then follow exactly the same line of argument as in (i) to see that
 \beaa
 u(t_n,x_n)-\cY_{t_n}^{t_n,x_n,\bar\dbP^n}\big(\tau_n, u_*(\tau_n,B^{t_n,x_n}_{\tau_n})\big)
 \le
 c_n - \eta',
 ~~\bar \dbP-\mbox{a.s.}
  \eeaa
where $\eta'$ is a positive constant independent of $n$. For large $n$, we have $c_n-\eta'<0$, and this is in contradiction with the dynamic programming principle.
\ep

\section{The dynamic programming principle}
\label{sect-dpp}
\setcounter{equation}{0}

In this section we prove Propositions \ref{prop-DP} and \ref{prop-u}.

\subsection{Regular conditional probability distributions}
\label{sect-rcpd}

The key tool to prove the dynamic programming principle is  the regular conditional probability distributions (r.c.p.d.), introduced by Stroock-Varadhan \cite{SV}. We adopt the notations of our
accompanying paper \cite{STZ09c}. For $0\le t \le s\le 1$, $\o \in \O^{t}$, $\tilde\o \in \O^{s}$, and $\cF^{t}_1-$measurable random variable $\xi$, define:
\bea
\label{xio}
\xi^{s,\o}(\tilde\o) := \xi(\o\otimes_{s}\tilde\o) ~\mbox{where}~ (\o\otimes_{s} \tilde\o) (r) := \o_r\1_{[t,s)}(r) + (\o_{s} + \tilde\o_r)\1_{[s, 1]}(r),~ r\in [t,1].
 \eea
 In particular, for any $\dbF^t-$stopping time $\t$, one can choose $s = \t(\o)$ and simplify the notation: $\o\otimes_\t \tilde\o :=\o\otimes_{\t(\o)}\tilde\o$. Clearly $\o\otimes_{\t}\tilde\o\in \O^{t}$ and,  for each $\o\in\O^t$, $\xi^{\t,\o}:=\xi^{\t(\o),\o}$ is $\cF^{\t(\o)}_1-$measurable. For each probability measure $\dbP$ on $(\O^{t}, \cF^{t}_1)$, by Stroock-Varadhan \cite{SV} there exist r.c.p.d. $\dbP^{\t,\o}$ for all $\o\in \O^{t}$ such that $\dbP^{\t,\o}$ is a probability measure on $(\O^{\t(\o)}, \cF^{\t(\o)}_1)$, and for all $\cF^{t}_1-$measurable $\dbP-$integrable random variable $\xi$:
\bea
\label{rcpd}
\dbE^\dbP[\xi|\cF^t_\t](\o) = \dbE^{\dbP^{\t,\o}}[\xi^{\t,\o}], &\mbox{for}& \dbP-\mbox{a.e.}~\o\in\O^{t}.
\eea
In particular, this implies that the mapping $\o\mapsto \dbE^{\dbP^{\t,\o}}[\xi^{\t,\o}]$ is $\cF^t_\t-$measurable.
Moreover, following the arguments in Lemmas 4.1 and 4.3 of \cite{STZ09c}, one can easily show that:

\begin{lem}
\label{lem-rcpd}
Let $t\in[0,1]$, $\t$ an $\dbF^t-$stopping time, and $\dbP\in \cP^{\k,t}_h$. Then:
\beaa
\mbox{for $\dbP-$a.e.}~\o\in\O^t:
&\dbP^{\t,\o}\in \cP^{\k,\t(\o)}_h ~\mbox{and}~ (\hat a^t)^{\t,\o}_r = \hat a^{\t(\o)}_r,&
dr\times d\dbP^{\t,\o} ~\mbox{on}~[\t(\o),1]\times \O^{\t(\o)}.
\eeaa
\end{lem}

\subsection{A weak partial dynamic programming principle}
\label{sect-weakDP1}

In this section, we prove the following result adapted from \cite{BT}.

\begin{lem}\label{lem-bt}
Under Assumptions \ref{assum-f} and \ref{assum-g}, for any $(t,x)$ and arbitrary $\dbF^t-$stopping times $\{\tau^\dbP,\dbP\in \cP_h^{\k,t}\}$:
 \beaa
 u(t,x)
 &\le&
 \sup_{\dbP \in \cP^{\k,t}_h} \cY^{t,x,\dbP}_t\big(\tau^\dbP, u^*(\tau^\dbP,B_{\tau^\dbP}^{t,x})\big).
 \eeaa
\end{lem}

\proof We shall prove the slightly stronger result:
\bea
\label{DP1-1}
\ba{c}
\cY^{t,x,\dbP}_t\big(1, g(B_1^{t,x})\big)
 \;\le\;
  \cY^{t,x,\dbP}_t\big(\tau^\dbP, \f(\tau^\dbP,B_{\tau^\dbP}^{t,x})\big)\\
  \mbox{for any}~\dbP\in \cP^{\k,t}_h~\mbox{and any Lebesgue measurable function}~ \f \ge u.
 \ea
 \eea
 Fix $\dbP$ and $\f$. For notation simplicity, we omit the dependence of $\tau^\dbP$ on $\dbP$. We first note that, by \ref{u<Lamda}, without loss of generality we may assume $|\f|\le \L$. Then Assumption \ref{assum-g} implies that $\cY^{t,x,\dbP}_t\big(\tau, \f(\tau,B_{\tau}^{t,x})\big)$ is well defined. By \reff{rcpd}, one can easily show that
 \beaa
 \cY_t^{t,x,\dbP}\big(1,g(B^{t,x}_1)\big)
 =
 \cY_t^{t,x,\dbP}\Big(\tau,\cY_\tau^{\tau(\o),B^{t,x}_{\tau}(\o),\dbP^{\t,\o}}\big(1,g(B^{\tau(\o),B^{t,x}_{\tau}(\o)}_1)\big)\Big)
 \eeaa
By Lemma \ref{lem-rcpd}, $\dbP^{\t,\o} \in \cP^{\k,\t(\o)}_h$, $\dbP-$a.e. $\o\in\O^t$. Then
\beaa
\cY_\tau^{\tau(\o),B^{t,x}_{\tau}(\o),\dbP^{\t,\o}}\big(1,g(B^{\tau(\o),B^{t,x}_{\tau}(\o)}_1)\big) \le u\big(\tau(\o),B^{t,x}_{\tau}(\o)\big) \le \f\big(\tau(\o),B^{t,x}_{\tau}(\o)\big),~\dbP-\mbox{a.e.}~\o\in\O^t.
\eeaa
It follows from the comparison result for BSDEs that
 \beaa
 \cY_t^{t,x,\dbP}\big(1,g(B^{t,x}_1)\big)
 &\le&
 \cY_t^{t,x,\dbP}\Big(\t, \f\big(\tau,B^{t,x}_\tau\big)\Big).
 \eeaa
This implies \reff{DP1-1}, and by the arbitrariness of $\dbP$, Lemma \ref{lem-bt} is proved.
\ep

\subsection{Concatenation of probability measures}

In preparation to the proof of Proposition \ref{prop-DP}, we introduce the concatenation of probability measures. For any $0\le t_0\le t\le 1$ and $\o\in \O^{t_0}$, denote $\o^t\in \O^t$ by $\o^t_s := \o_s-\o_t$, $s\in [t,1]$. For any $\dbP_1 = \dbP^{\a^1}\in \cP^{\k,t_0}_h$, $\dbP_2 = \dbP^{\a^2}\in \cP^{\k,t}_h$, let $\dbP:= \dbP_1\otimes_t \dbP_2$ denote the probability measure $\dbP^\a$, where
\beaa
\a_s (\o) := \a^1_s(\o) \1_{[t_0,t]}(s) + \a^2_s(\o^t)\1_{[t,1]}(s), ~~\o\in\O^{t_0}.
\eeaa

\begin{lem}
\label{lem-concatenation}
Let $\dbP  := \dbP_1\otimes_t \dbP_2$ be as defined above. Then, under Assumption \ref{assum-f},
\bea
\label{concatenation}
\dbP\in \cP^{\k,t_0}_h,~~
\dbP = \dbP_1 ~~\mbox{on}~~\cF^{t_0}_t, &\mbox{and}& \dbP^{t,\o} = \dbP_2 ~~\mbox{for}~~\dbP_1-\mbox{a.e.}~\o\in \O^{t_0}.
\eea
\end{lem}

\proof First by \reff{ellipticity-t}, we have $\underline a_{\dbP_i} \le \a^i \le \overline a_{\dbP_i}$, $i=1,2$. Then $\underline a_{\dbP_1} \wedge \underline a_{\dbP_2} \le \a \le \overline a_{\dbP_1}\vee \overline a_{\dbP_1}$. In particular, this implies that $\int_{t_0}^1 |\a_s|ds <\infty$. Then $\dbP\in \overline\cP^{t_0}_S$ and $\underline a_{\dbP_1} \wedge \underline a_{\dbP_2} \le \hat a \le \overline a_{\dbP_1}\vee \overline a_{\dbP_2}$, $\dbP-$a.s. The two last claims in \reff{concatenation} are obvious, and imply that:
$$
\begin{array}{rcl}
 \dbE^\dbP\Big[\big(\int_{t_0}^1|\hat f^{t_0,0}_s|^\k ds\big)^{2/\k}\Big]
 &\le&
 C_\k\dbE^\dbP\Big[\big(\int_{t_0}^t|\hat f^{t_0,0}_s|^\k ds\big)^{2/\k}                    + \big(\int_t^1 |\hat f^{t_0,0}_s|^\k ds\big)^{2/\k}
              \Big]
 \\
 &=&
 C_\k\left(\dbE^{\dbP_1}\Big[\big(\int_{t_0}^t|\hat f^{t_0,0}_s|^\k ds
                             \big)^{2/\k}
                        \Big]
            + \dbE^{\dbP_1}
              \Big[\dbE^{\dbP_2}\big[\big(\int_t^1|\hat f^{t,0}_s|^\k ds
                                     \big)^{2/\k}
                                \big]
              \Big]
       \right)
 \\
 &=&
 C_\k\left(\dbE^{\dbP_1}\Big[\big(\int_{t_0}^t|\hat f^{t_0,0}_s|^\k ds
                             \big)^{2/\k}
                        \Big]
           + \dbE^{\dbP_2}\Big[\big(\int_t^1 |\hat f^{t,0}_s|^\k ds
                               \big)^{2/\k}
                          \Big]
     \right) < \infty.
\end{array}
$$
This implies that  $\dbP\in \cP^{\k,t_0}_h$.
\ep

\subsection{Dynamic programming and regularity}
\label{subsect-DPlsc}

We first prove the dynamic programming principle of Proposition \ref{prop-DP} for stopping times taking countably many values. From this, we will deduce the lower-semicontinuity of $u$ stated in Proposition \ref{prop-u}, which in turn provides Proposition \ref{prop-DP} by passing to limits. 

\begin{lem}
\label{lem-strongDP}
Proposition \ref{prop-DP} holds true under the additional condition that each $\t^\dbP$ takes countable many values.
\end{lem}

\proof (i) We first observe that the lower semicontinuity of $g$ implies that
 \bea\label{Y lsc in x}
 x \longmapsto \cY^{t,x,\dbP}_{t}(1, g(B^{t,x}_1))
 &\mbox{is lower-semicontinuous for all}&
 \dbP\in\cP^{\k,t}_h.
 \eea
This is a direct consequence of the stability and comparison principle of BSDEs. Then, for all fixed $(t,x)$, and all sequence $(x_n)_{n\ge 1}$ converging to $x$, it follows that:
 \beaa
 u(t,x)
 =\sup_{\dbP\in\cP^{\k,t}_h}\cY^{t,x,\dbP}_t\big(1,g(X^{t,x}_1)\big)
 \le
 \sup_{\dbP\in\cP^{\k,t}_h}\liminf_{n\to\infty}\cY^{t,x_n,\dbP}_t\big(1,g(X^{t,x_n}_1)\big)
 \le
 \liminf_{n\to\infty}u(t,x_n).
 \eeaa
Hence $u(t,.)$ is lower-semicontinuous, and therefore measurable.
\\
(ii) We now fix $(t_0,x_0)$ and prove the result at this point. Let $\t$ be an $\dbF^{t_0}-$stopping time with values in $\{t_k, k\ge 1\}\subset [t_0,1]$. Since $u(t_k,.)$ is measurable, we deduce that $u(\t, B_\t^{t_0,x_0}) = \sum_{k\ge 1} u(t_k, B^{t_0,x_0}_{t_k})\1_{\{\t=t_k\}}$ is $\cF_\t-$measurable. Then, it follows from \reff{DP1-1} that
\beaa
u(t_0,x_0) &\le& \sup_{\dbP \in \cP^{\k,t_0}_h} \cY^{t_0,x_0,\dbP}_{t_0}\big(\tau^\dbP, u(\tau^\dbP,B_{\tau^\dbP}^{t_0,x_0})\big).
\eeaa
(iii) To complete the proof, we fix $\dbP\in \cP^{\k, t_0}_h$, denote $\t:=\t^\dbP$, and proceed in four steps to show that
\bea
\label{DPstrong1-1}
  \cY^{t_0,x_0,\dbP}_{t_0}\big(\tau, u(\tau,B_{\tau}^{t_0,x_0})\big) &\le & u(t_0,x_0).
\eea
{\it Step 1.} We first fix  $t\in (t_0,1]$, and show that,
\bea
\label{DPstrong1-2}
\cY^{t_0,x_0,\dbP}_{t_0}(t, \f(B^{t_0,x_0}_{t})) \le u(t_0,x_0),
\eea
for any continuous function $\varphi:\dbR^d\longrightarrow\dbR$ such that  $-\L(t,\cd) \le \f(\cd) \le u(t,\cd)$.
Indeed, for any $\dbP^t\in \cP^{\k,t}_h$, by the lower-semicontinuity property \reff{Y lsc in x}, we may argue exactly as in Step 2 of the proof of Theorem 3.1 in \cite{BT} to deduce that, for every $\eps>0$,
 there exist sequences $(x_i,r_i)_{i\ge 1}\subset \dbR^d\times(0,1]$ and $\dbP_i\in\cP_h^{\k,t}$, $i\ge 1$ such that
 \beaa
 \cY^{t,\cd,\dbP_i}_{t}(1, g(B^{t,\cd}_1))\ge \f(t,\cd)-\eps
 ~\mbox{on}~
 Q_i:=\{x'\in\dbR^d: |x'-x_i|<r_i\},
 &\mbox{and}&
 \cup_{i\ge 1}Q_i = \dbR^d.
 \eeaa
This provides a disjoint partition $(A_i)_{i\ge 1}$ of $\dbR^d$ defined by $A_{i}:=Q_{i}\setminus \cup_{j<i} Q_j$. Set
\beaa
E_i := \{B_{t}^{t_0,x_0}\in A_i\}, ~~i\ge 1, &\mbox{and}& \bar E_n := \cup_{i> n} E_i, ~~n\ge 1.
\eeaa
Then $\{E_i, 1\le i\le n\}$ and $\bar E^n$ form a partition of $\O$ and $\lim_{n\to\infty}\dbP(\bar E_n)=0$.
Define
\bea
\label{barPndef}
\bar\dbP^n(E) := \sum_{i=1}^n (\dbP\otimes_t\dbP_i)(E\cap E_i) + \dbP(E\cap \bar E_n)  ~~\mbox{for all}~~E\in \cF^{t_0}_1.
\eea
Combining the arguments for  \reff{barPN-in-cP} and Lemma \ref{lem-concatenation}, one can easily show that
\bea
\label{barPn}
\bar\dbP^n\in \cP^{\k,t_0}_h(t,\dbP) &\mbox{and}& (\bar\dbP^n)^{t,\o} = \dbP_i,~~\dbP-\mbox{a.e.}~\o\in E_i,~~1\le i\le n.
\eea
This implies that, for $1\le i\le n$ and $\dbP-$a.e. $\o\in E_i$,
\beaa
\cY^{t_0,x_0, \bar\dbP^n}_{t}(1, g(B^{t_0,x_0}_1))(\o)
=\cY^{t,B^{t_0,x_0}_t(\o),\dbP_i}_{t}(1, g(B^{t,B^{t_0,x_0}(\o)_t}_1))
   \ge \f\big(B^{t_0,x_0}_t(\o)\big)-\e,
\eeaa
and, by the comparison result for BSDEs:
 \beaa
 u(t_0,x_0)
 &\ge& \cY^{t_0,x_0, \bar\dbP^n}_{t_0}(1, g(B^{t_0,x_0}_1)) = \cY^{t_0,x_0, \dbP}_{t_0}\Big(t,\cY^{t_0,x_0, \bar\dbP^n}_{t}(1, g(B^{t_0,x_0}_1))\Big)
 \\
 &\ge&
 \cY^{t_0,x_0, \dbP}_{t_0}\Big(t,\big(\f(B^{t_0,x_0}_t)-\e\big)\1_{(\bar E_n)^c} + \cY^{t_0,x_0, \bar\dbP^n}_{t}(1, g(B^{t_0,x_0}_1))\1_{\bar E_n}\Big).
 \eeaa
By the stability of BSDEs and the arbitrariness of $\eps>0$, this proves \reff{DPstrong1-2}.

\noindent {\it Step 2.} Since $u(t,\cd)$ is lower semi-continuous, there exist continuous functions $\{\f_n, n\ge 1\}$ such that $\f_n\uparrow u(t,\cd)$. Without loss of generality we may assume $\f_n \ge -\L$. Since \reff{DPstrong1-2} holds for each $\f_n$, we obtain \reff{DPstrong1-1} for $\t=t$ by monotone convergence.

\noindent {\it Step 3.} Assume $\t$ takes finitely many values $t_0<t_1<\cds<t_n \le 1$. Note that, $\dbP-$a.s.
\beaa
&&\cY^{t_0,x_0,\dbP}_{\t\wedge t_{n-1}}\big(\t, u(\tau,B_{\tau}^{t_0,x_0})\big)\\
 &=& \cY^{t_0,x_0,\dbP}_{\t}\big(\t, u(\tau,B_{\tau}^{t_0,x_0})\big)\1_{\{\t\le t_{n-1}\}} + \cY^{t_0,x_0,\dbP}_{t_{n-1}}\big(\t, u(\tau,B_{\tau}^{t_0,x_0})\big)\1_{\{\t > t_{n-1}\}}\\
 &=& u(\tau,B_{\tau}^{t_0,x_0})\1_{\{\t\le t_{n-1}\}} + \cY^{t_{n-1},B^{t_0,x_0}_{t_{n-1}}(\o),(\dbP)^{t_{n-1},\o}}_{t_{n-1}}\big(t_n, u(t_n,B_{t_n}^{t_{n-1},B^{t_0,x_0}_{t_{n-1}}(\o)})\big)\1_{\{\t > t_{n-1}\}}
\eeaa
By Lemma \ref{lem-rcpd}, $(\dbP)^{t_{n-1},\o}\in \cP^{\k,t_{n-1}}_h$, $\dbP-$a.s. Then by Step 2 we have
\beaa
\cY^{t_0,x_0,\dbP}_{\t\wedge t_{n-1}}\big(\t, u(\tau,B_{\tau}^{t_0,x_0})\big)&\le& u(\tau,B_{\tau}^{t_0,x_0})\1_{\{\t\le t_{n-1}\}} + u(t_{n-1},B^{t_0,x_0}_{t_{n-1}})\1_{\{\t > t_{n-1}\}}\\
 &=& u(\t\wedge t_{n-1}, B^{t_0,x_0}_{\t\wedge t_{n-1}}).
\eeaa
Then, by the comparison principle of BSDE,
\beaa
\cY^{t_0,x_0,\dbP}_{t_0}\big(\tau, u(\tau,B_{\tau}^{t_0,x_0})\big) &=& \cY^{t_0,x_0,\dbP}_{t_0}\Big(\t\wedge t_{n-1}, \cY^{t_0,x_0,\dbP}_{\t\wedge t_{n-1}}\big(\t, u(\tau,B_{\tau}^{t_0,x_0})\big)\Big)\\
&\le& \cY^{t_0,x_0,\dbP}_{t_0}\Big(\t\wedge t_{n-1}, u(\t\wedge t_{n-1}, B^{t_0,x_0}_{\t\wedge t_{n-1}})\Big).
\eeaa
Continuing this backward induction provides \reff{DPstrong1-1}.

\noindent {\it Step 4.} Now assume $\t$ takes countable many values $\{t_k, k\ge 1\}$. Denote $\t_n := \sum_{k=1}^n t_k \1_{\{\t=t_k\}} + \1_{\{\t \neq t_k, 1\le k\le n\}}$. Clearly $\t_n$ is still an $\dbF^{t_0}-$stopping time. By Step 3,
\beaa
\cY^{t_0,x_0,\dbP}_{t_0}\big(\tau_n, u(\tau_n,B_{\tau_n}^{t_0,x_0})\big) &\le & u(t_0,x_0).
\eeaa
For each $\o\in \O^{t_0}$, we have  $\t_n(\o) = \t(\o)$, for sufficiently large $n$. Then $u(\tau_n(\o),B_{\tau_n}^{t_0,x_0}(\o)) = u(\tau(\o),B_{\tau}^{t_0,x_0}(\o))$, and \reff{DPstrong1-1} follows from the stability of BSDEs.
\ep

\vspace{5mm}

As a consequence of Lemma \ref{lem-strongDP}, we can now prove that $u$ is lower-semicontinuous.

\vspace{5mm}

\noindent {\it Proof of Proposition \ref{prop-u}.}\quad Recall the $\cY^\dbP(\t,\xi)$ defined in \reff{BSDEP}, and define
\bea
\label{J}
J(t,x,\dbP) := \dbE^\dbP\Big[\cY^{\dbP}_t(1, g(x + B_1-B_t))\Big] &\mbox{for all}& t, x, ~\mbox{and}~ \dbP\in \cP^{\k}_h.
\eea
(i) We first prove that
\bea
\label{u=J}
u(t,x) = \sup_{\dbP\in \cP^{\k}_h} J(t,x,\dbP).
\eea
To see this, we first observe that, for any $\dbP\in \cP^{\k}_h$, it follows from Lemma \ref{lem-rcpd} that
\beaa
\cY^{\dbP}_t(1, g(x + B_1-B_t))(\o) = \cY^{t,x, \dbP^{t,\o}}_t(1, g(B^{t,x}_1)) \le u(t,x) &\mbox{for}& \dbP-\mbox{a.e.}~\o\in\O.
\eeaa
Then $J(t,x,\dbP)\le u(t,x)$ for any $\dbP\in \cP^{\k}_h$.

On the other hand, for any $\dbP_2\in \cP^{\k,t}_h$, choose arbitrary $\dbP_1\in\cP^{\k}_h$ and let $\dbP := \dbP_1 \otimes_t \dbP_2$. Then $\dbP\in \cP^{\k}_h$ and, by \reff{concatenation},
\beaa
\cY^{\dbP}_t(1, g(x + B_1-B_t))(\o) =\cY^{t,x, \dbP^{t,\o}}_t(1, g(B^{t,x}_1)) = \cY^{t,x, \dbP_2}_t(1, g(B^{t,x}_1)) &\mbox{for}& \dbP-\mbox{a.e.}~\o\in\O.
\eeaa
This implies that $J(t,x,\dbP) = \cY^{t,x, \dbP_2}_t(1, g(B^{t,x}_1))$ and thus
$
u(t,x) \le \sup_{\dbP\in \cP^{\k}_h} J(t,x,\dbP).
$

\ms

\noindent (ii) We now prove that the lower-semicontinuity of $g$ implies that:
 \bea
 \label{Jlsc}
 (t,x)\longmapsto J(t,x,\dbP)
 &\mbox{is lower-semicontinuous for any}&
 \dbP\in \cP^{\k}_h.
 \eea
which obviously implies the lower-semicontinuity of $u$ in view of \reff{J}.

For $(t,x)\in[0,1]\times\dbR^d$ and $\dbP\in\cP^{\k}_h$, let $(t_n,x_n)_{n\ge 1}$ be a sequence in $[0,1]\times\dbR^d$ such that
$
 (t_n,x_n)\longrightarrow (t,x)$.
Denote, for each $n$,
\beaa
\xi_n := \inf_{k\ge n} g( x_k+ B_1-B_{t_k}),&& f^n_s(y,z) := \inf_{k\ge n} f(s, x_k+B_s-B_{t_k}, y,z, \hat a_s),\\
\xi_\infty := \lim_{n\to\infty} \xi_n, && f^\infty := \lim_{n\to\infty} f^n,
\eeaa
and, for $1\le n \le \infty$, let $(\cY^n, \cZ^n)$ denote the solution to the following BSDE:
\beaa
\cY^n_s = \xi_n - \int_s^{1} f^n_r(\cY^n_r, \cZ^n_r) dr - \int_s^{1} \cZ^n_r dB_r,~~ t\le s\le 1, ~~\dbP-\mbox{a.s.}
\eeaa
By Assumptions \ref{assum-f} and \ref{assum-g}, $g$ and the modulus of continuity $\rho$ of $f$ have polynomial growth in $x$. Then  there exist some constants  $C$ and $p$ such that
\bea
\label{polynomial}
\sup_{n\ge 1}\Big\{|\xi_n| + |f^n_t(0,0)|\Big\}\le |\hat f^{0,0}_r| + C\Big(\sup_{k\ge 1}|x_k|^p + \sup_{0\le t\le 1}|B_t|^p\Big).
\eea
Moreover, $\hat a$ has upper bound $\overline a_\dbP$, $\dbP-$a.s. then it follows from the Lipschitz conditions of $f$ that  the above BSDE has a unique solution for each $n$, and
\beaa
\lim_{n\to\infty}\dbE^\dbP[\cY^n_t]= \dbE^\dbP[\cY^\infty_t].
 \eeaa
 By the lower semi-continuity of $g$ and the uniform continuity of $f$ in $x$ in \reff{fLip}, we have $\xi_\infty \ge g(x+B_1-B_t)$ and $f^\infty_s(y,z) = f(s, x+B_1-B_s,y,z, \hat a_s)$, $\dbP-$a.s.  Then by  the comparison principle of BSDEs  one can easily see that
\beaa
\liminf_{n\to\infty}J(t_n,x_n,\dbP) \ge \lim_{n\to\infty}\dbE^\dbP[\cY^n_t]= \dbE^\dbP[\cY^\infty_t] \ge \dbE^\dbP[\cY^{\dbP}_t\big(1, g(x+B_1-B_t)\big)\Big]=J(t,x,\dbP).
\eeaa
This proves the lower-semicontinuity of $J$ for any fixed $\dbP\in \cP^{\k}_h$.
\ep

\vspace{5mm}

We now can prove the dynamic programming principle for arbitrary stopping times.

\vspace{5mm}

\noindent {\it Proof of Proposition \ref{prop-DP}}\quad
For any $(t,x)$, $\dbP\in \cP^{\k,t}_h$, $\dbF^t-$stopping time $\t$, and any $n$, denote
\beaa
\t_n := \sum_{i=1}^n \frac{i}{n}\1_{[i-\frac{1}{n}, \frac{i}{n})}(\t) + \1_{\{\t=1\}}.
\eeaa
Then $\t_n$ is an $\dbF^t-$stopping time, $\t_n \ge \t$, and $\t_n \to \t$. By  Lemma \ref{lem-strongDP}, together with Proposition \ref{prop-u}, we have
\beaa
\cY^{t,x,\dbP}_{t}\big(\tau_n, u(\tau_n,B_{\tau_n}^{t,x})\big) &\le & u(t,x)
\eeaa
Since $u$ is lower-semicontinuous, $\liminf_{n\to\infty} u(\tau_n,B_{\tau_n}^{t,x}) \ge u(\t, B^{t,x}_\t)$. Then it follows from the comparison and the stability of BSDEs that
 \beaa
 u(t,x)
 &\ge&
  \cY^{t,x,\dbP}_t\left(\t,u(\t,B^{t,x}_\t)\right).
 \eeaa
Finally, $u$ is measurable since it is lower-semicontinuous. Then \reff{DP1-1} provides the opposite inequality.
\ep

\section{Appendix}
\label{sect-appendix}
\setcounter{equation}{0}

\subsection{Non-uniqueness in $\dbL^2(\dbP_0)$ of the 2BSDE \reff{2BSDEc}}
\label{sect-linear 2BSDE}

In this section, we provide an example which shows the importance of the constraints imposed in \cite{CSTV} to obtain uniqueness.

\begin{eg}
\label{eg-counterexample}
Consider the following 2 dimensional forward SDEs:
 \bea
 \label{counterexample}
 \left\{\ba{lll}
 \dis Y_t = -\int_0^t \frac{3Y_s}{1-s}ds + \int_0^t \frac{X_s}{\sqrt{1-s}}dB_s,\\
 \dis X_t = 1-\int_0^t \frac{3(1+c^2)X_s}{2c^2(1-s)}ds + \int_0^t \frac{3Y_s}{c\sqrt{1-s}}dB_s,
 \ea\right. \dbP_0-\mbox{a.s.}
 \eea
Clearly, (\ref{counterexample}) is well-posed on $[0,1)$. Denote
 \beaa
 Z_t := \frac{X_t}{\sqrt{1-t}};\q \G_t := \frac{3Y_t}{c(1-t)};\q A_t := -\big(\frac{3}{2c^2}+1\big)\frac{X_t}{(1-t)^{3/2}}.
 \eeaa
Then $(Y,Z,\G,A)$ is a nonzero solution to 2BSDE (\ref{2BSDEc}).
\end{eg}

\proof
First, applying It\^o's formula one can check straightforwardly that $(Y,Z,\G,A)$ satisfies the SDEs in (\ref{2BSDEc}).
Notice that
 \beaa
 R_t := \frac{3}{c^2}Y_t^2 + X_t^2
 &\mbox{satisfies}&
 dR_t = -\frac{3R_t}{1-t} dt + (\cds)dB_t,
 \eeaa
by It\^o's formula. Since $R_0=1$,
 \beaa
 \dbE^{\dbP_0}[R_t] = 1 -3 \int_0^t \frac{\dbE^{\dbP_0}[R_s]}{1-s}ds
 &\mbox{and thus}
 & \dbE^{\dbP_0}[R_t] = (1-t)^3,~~\mbox{for all}~0\le t<1.
 \eeaa
Then one can easily see that,
 \beaa
 \sup_{0\le t<1} \dbE^{\dbP_0}\big[|\G_t|^2 + |A_t|^2\big] \le C\dbE^{\dbP_0}\Big[\frac{|Y_t|^2}{(1-t)^2} + \frac{|X_t|^2}{(1-t)^3}\Big] \le C,
 \eeaa
which, together with \reff{2BSDEc}, also implies that
 \beaa
 \dbE^{\dbP_0}\Big[\sup_{0\le t<1} [|Y_t|^2 + |Z_t|^2\Big] \le C.
 \eeaa
Finally, we prove that
 \bea
 \label{Y1}
 \lim_{t\uparrow 1} Y_t = 0, ~\dbP_0-\mbox{a.s.}
 \eea
In fact, for any $t<T<1$, by Burkholder-Davis-Gundy inequality we have
 \beaa
 \dbE\Big[\sup_{t\le s\le T}|Y_s|^2\Big]
 &\le&
 C\dbE\Big[|Y_t|^2 + \int_t^T \frac{|Y_s|^2}{(1-s)^2}ds
                                                           + \int_t^T \frac{|X_s|^2}{1-s}ds\Big]\\
 &\le& C\Big((1-t)^3 + \int_t^T \big((1-s)+(1-s)^2\big)ds\Big)\le
 C(1-t)^2.
 \eeaa
Let $T\uparrow 1$ and apply the monotone convergence Theorem, we get
 \beaa
 \dbE\Big[\sup_{t\le s< 1}|Y_s|^2\Big]\le C(1-t)^2.
 \eeaa
Then $\sup_{t\le s< 1}|Y_s|^2\downarrow 0$, as $t\uparrow 1$, $\dbP_0-$a.s. by the decrease of $\sup_{t\le s< 1}|Y_s|^2$ in $t$, and we deduce \reff{Y1}.
 \ep

\subsection{Proof of Lemma \ref{lem-BSDEest}}
\label{proof-lem-BSDEest}

If the a priori estimates \reff{BSDEYest} and \reff{BSDEZest} hold, then by the martingale representation property \reff{0-1 MRP}, the Lipschitz conditions \reff{FLip-yz}, and the integrability assumption of $\hat F^0$ in \reff{ellipticity}, following the standard arguments one can easily show that BSDE \reff{BSDEP} has a unique solution.

We now prove \reff{BSDEYest} and \reff{BSDEZest}. For notational simplicity in the proof we drop the superscripts $\dbP$ in $(\cY^\dbP, \cZ^\dbP)$.
By the Lipschitz conditions \reff{FLip-yz}, there exist bounded processes $\l, \eta$ such that
\bea
\label{cYt}
\cY_t = \xi + \int_t^1 \big(\hat F^0_s+\l_s \cY_s + \eta_s \hat a_s^{1\slash 2} \cZ_s \big) ds
         - \int_t^1 \cZ_s dB_s, ~0\le t\le 1, ~\dbP-\mbox{a.s.}
\eea
Define $M$ by \reff{M}. By It\^o's formula, we have:
\bea
\label{McY}
d\big(M_t \cY_t\big)
= -M_t\hat F^0_t dt +
M_t\big( \cZ_t - \cY_t \eta_t{\hat a}^{-1/2}_t \big) dB_t,
~~0\le t\le 1, ~\dbP-\mbox{a.s.}
\eea
Then, using standard localization arguments if necessary:
\beaa
\cY_{t}
= M_t^{-1}\dbE^\dbP_t\Big[M_1\xi + \int_t^1 M_s \hat F^0_s ds\Big],~~ 0\le t\le 1, ~\dbP-\mbox{a.s.}
\eeaa
It follows from \reff{Mest} that, for $1<\k \le 2$,
\beaa
|\cY_{t}|
&\le& \dbE^\dbP_t\Big[\sup_{t\le s\le 1}(M_t^{-1}M_s) \big(|\xi| + \int_t^1 |\hat F^0_s| ds\big)\Big]\\
&\le& C_\k \Big(\dbE^\dbP_t\Big[|\xi|^\k + \int_t^1 |\hat F^0_s|^\k ds\big)\Big]\Big)^{1/\k},~~ 0\le t\le 1, ~\dbP-\mbox{a.s.}
\eeaa
This proves \reff{BSDEYest}.

Finally, applying It\^{o}'s formula on $\cY_t^2$ and following standard arguments we have
\beaa
&&\dbE^\dbP\Big[\int_0^1 |\hat a_t^{1/2} \cZ_t|^2 dt\Big] \le C\dbE^\dbP\Big[|\xi|^2 + \sup_{0\le t\le 1}|\cY_t| \int_0^1 |\hat F^0_t|dt\Big]\\
&\le& C\dbE^\dbP\Big[\sup_{0\le t\le 1}|\cY_t|^2 + \Big(\int_0^1 |\hat F^0_t|dt\Big)^2\Big]\le C\dbE^\dbP\Big[\sup_{0\le t\le 1}|\cY_t|^2 + \Big(\int_0^1 |\hat F^0_t|^\k dt\Big)^{2/\k}\Big].
\eeaa
This, combing with \reff{BSDEYest}, proves \reff{BSDEZest}.
\ep

\subsection{Proof of \reff{barPN-in-cP}}
\label{proof-barPN-in-cP}

By the definition of ${\cP_H^\k}$, we have $\dbP = \dbP^\a$, $\dbP'_1 = \dbP^{\a^1}$, and $\dbP'_2 = \dbP^{\a^2}$ for $\dbF-$progressively measurable processes $\a, \a^1,\a^2$ taking values in $\dbS^{>0}_d$. Since $\dbP, \dbP'_1,\dbP'_2\in{\cP_H^\k}$, by \reff{ellipticity} there exist $\underline\a, \overline\a, \underline\a^i, \overline\a^i\in \dbS^{>0}_d$ such that
\beaa
\underline\a \le \a \le \overline\a, ~~~\underline\a^i \le \a^i \le \overline\a^i,~~dt\times d\dbP_0-\mbox{a.s.}
\eeaa
Since $\dbP'_i\in{\cP_H^\k}(t,\dbP)$, it is clear that $\a = \a^i$, $ds\times d\dbP_0-$a.s. on $[0,t]\times\O$. Recall \reff{Xa}, then
\beaa
\a^*_s(\o) := \a_s(\o) \1_{[0,t)}(s) + \Big(\a^1_s(\o) \1_{\{X^\a \in E_1\}}(\o)+ \a^2_s(\o) \1_{\{X^\a \in E_2\}}(\o)\Big)\1_{[t,1]}(s),
~s\in[0,1],
\eeaa
is $\dbF-$progressively measurable and satisfies:
 \beaa
 0< \underline\a \wedge \underline\a^1\wedge \underline\a^2
 \le
 \a^*
 \le
 \overline\a \vee \overline\a^1\vee \overline\a^2.
 \eeaa
Following a line by line analogy of the proof of Claim 4.19 in \cite{STZ09c}, which in turn uses the arguments in the proof of Lemma 4.1  in \cite{STZ09c}, we see that $\dbP'=\dbP^{\a^*} \in \overline\cP_S$.  Moreover,
 \beaa
 \dbE^{\dbP'}\Big[\int_0^1 |\hat F^0_s|^2 ds \Big]
 &=&
 \dbE^{\dbP}\Big[\int_0^t |\hat F^0_s|^2 ds \Big]
 + \dbE^{\dbP'_1}\Big[\int_t^1 |\hat F^0_s|^2 ds \1_{E_1}\Big]
 + \dbE^{\dbP'_2}\Big[\int_t^1 |\hat F^0_s|^2 ds \1_{E_2}\Big]\\
 &\le&
 \dbE^{\dbP}\Big[\int_0^1 |\hat F^0_s|^2 ds \Big]
 + \dbE^{\dbP'_1}\Big[\int_0^1 |\hat F^0_s|^2 ds \1_{E_1}\Big]
 + \dbE^{\dbP'_2}\Big[\int_0^1 |\hat F^0_s|^2 ds \1_{E_2}\Big]
 <\infty.
 \eeaa
Then $\dbP'\in {\cP_H^\k}$. Obviously, $\dbP' = \dbP$ on $\cF_t$. This proves that $\dbP'\in{\cP_H^\k}(t,\dbP)$.
\ep

\end{document}